\documentclass[11pt]{amsart}

\usepackage{amsthm,amssymb,amsmath,calc,eucal,ifthen}
\usepackage[matrix,arrow]{xy}

\newtheorem{theorem}{Theorem}[section]
\newtheorem{lemma}[theorem]{Lemma}
\newtheorem{proposition}[theorem]{Proposition}
\newtheorem{corollary}[theorem]{Corollary}

\newtheorem*{theorem*}{Theorem}
\newtheorem*{proposition*}{Proposition}

\theoremstyle{remark}
\newtheorem{remark}[theorem]{Remark}

\newtheoremstyle{rmdefinition}{}{}{\upshape}{}{\bfseries}{.}{ }{}

\theoremstyle{rmdefinition}
\newtheorem{definition}[theorem]{Definition}

\newcommand{\OMIT}[1]{}

\newcommand{\abs}[1]{\lvert #1\rvert}

\newcommand{\derived}[2][1]{\ifthenelse{\equal{#1}{1}}{{#2}^{\prime}}{%
\ifthenelse{\equal{#1}{2}}{{#2}^{\prime\prime}}{\ifthenelse{\equal{#1}{3 
  }}
{{#2}^{\prime\prime\prime}}{{#2}^{(#1)}}}}}
\newcommand{\subgroup}{<}

\newcommand{\norm}[1]{\|#1\|}

\newcommand{\calA}{\mathcal{A}}
\newcommand{\calC}{\mathcal{C}}
\newcommand{\calD}{\mathcal{D}}
\newcommand{\calG}{\mathcal{G}}
\newcommand{\calH}{\mathcal{H}}
\newcommand{\calK}{\mathcal{K}}

\newcommand{\Fsupplus}{{\overline{F}}^+\!}
\newcommand{\Finfplus}{{\underline{F}}^+\!}

\newcommand{\fsup}{\overline{f}}
\newcommand{\finf}{\underline{f}}
\newcommand{\fsupplus}{{\overline{f}}^+\!}
\newcommand{\finfplus}{{\underline{f}}^+\!}
\newcommand{\fplus}{f^+\!}

\DeclareMathOperator{\Aut}{Aut}

\DeclareMathOperator{\im}{im}
\DeclareMathOperator{\pr}{pr}
\DeclareMathOperator{\tr}{tr}
\DeclareMathOperator{\fin}{fin}
\DeclareMathOperator{\spec}{spec}
\DeclareMathOperator{\specpoint}{spec_{\text{\rm point}}}
\DeclareMathOperator{\Tr}{Tr}

\newcommand{\hyp}{{\mathbb H}}

\newcommand{\R}{{\mathbb R}}
\newcommand{\C}{{\mathbb C}}
\newcommand{\Z}{{\mathbb Z}}
\newcommand{\Qbar}{\overline{\mathbb{Q}}}

\renewcommand{\th}{\text{th}}

\begin{document}

\title[On eigenvalues of generalized Harper operators on graphs]
{Arithmetic properties of eigenvalues of generalized Harper operators
   on graphs}

\author{J\'ozef Dodziuk}
\address{Ph.D.\ Program in Mathematics, Graduate Center of CUNY,  
New~York,
   NY 10016, USA.}
\email{jozek@derham.math.qc.edu}

\author{Varghese Mathai}
\address{Department of Mathematics, University of Adelaide, Adelaide  
5005,
   Australia}
\email{vmathai@maths.adelaide.edu.au}

\author{Stuart Yates}
\address{Department of Mathematics, University of Adelaide, Adelaide  
5005,
   Australia and Max Planck Institut f\"ur Mathematik, Bonn, Germany}
\email{syates@maths.adelaide.edu.au}


\thanks{%
   The second and third authors acknowledge support from the Australian
   Research Council.%
}

\begin{abstract}
   Let $\Qbar$ denote the field of algebraic numbers in $\C$.
   A discrete group $G$ is said to have the $\sigma$-multiplier algebraic
   eigenvalue property, if for every matrix $A\in M_d(\Qbar(G,\sigma))$,
   regarded as an operator on $l^2(G)^d$, the eigenvalues of $A$ are
   algebraic numbers, where $\sigma \in Z^2(G,\mathcal U(\Qbar))$
   is an algebraic multiplier, and $\mathcal U(\Qbar)$
   denotes the unitary elements of $\Qbar$.  Such operators
   include the Harper operator and the discrete magnetic Laplacian that
   occur in solid state physics. We prove that any finitely generated
   amenable, free or surface group has this property for any 
   algebraic multiplier $\sigma$.
   In the special
   case when $\sigma$ is rational ($\sigma^n$=1 for some positive
   integer $n$) this property holds for a larger class of groups
   $\mathcal K$ containing free groups and amenable groups, and closed
   under taking directed unions and extensions with amenable quotients.
   Included in the paper are proofs of other spectral properties of
   such operators.
\end{abstract}

\maketitle

\section{Introduction}

This paper is concerned with number theoretic properties of
eigenvalues of self adjoint matrix operators
that are associated with weight functions on a graph equipped
with a free action of a discrete group. These operators form
generalizations of the Harper operator and the discrete magnetic
Laplacian (DML) on such graphs, as defined by Sunada in \cite{Sun}.

The Harper operator and DML over the Cayley graph of ${\mathbb Z}^2$
arise as the Hamiltonian in discrete models of the behaviour of free
electrons in the presence of a magnetic field, where the strength of
the magnetic field is encoded in the weight function.  When the weight
function is trivial, the Harper operator and the DML reduce to the
Random Walk operator and the discrete Laplacian respectively.  The DML
is in particular the Hamiltonian in a discrete model of the integer
quantum Hall effect (see for example \cite{Bel}); when the graph is
the Cayley graph of a cocompact Fuchsian group, the DML becomes the
Hamiltonian in a discrete model of the fractional quantum Hall effect
(\cite{CHMM98}, \cite{CHM} and \cite{MM01}).  It has also been studied
in the context of noncommutative Bloch theory (\cite{MM99}, \cite{Sh}.)

The Harper operator and DML can be thought of as particular examples
of weighted sums of twisted right translations by elements of the
group; alternatively, they can be regarded as matrices over the group
algebra twisted by a $2$-cocycle, acting by (twisted) left
multiplication. As such, in section \ref{sec:algmult} we
generalize results of \cite{dlmsy} to demonstrate that
such operators associated with algebraic weight functions have only
algebraic eigenvalues whenever the group is in a class of groups
containing all free groups, finitely generated amenable groups and
fundamental groups of closed Riemann surfaces.

When the multiplier associated to the weight function is rational, the
algebraicity of eigenvalues extends to groups in a larger class
$\mathcal K$, defined in section \ref{sec:cyclic}. The class $\mathcal  
K$
contains all free groups, discrete amenable groups, and groups in
the Linnell class $\mathcal C$; in particular it includes cocompact
Fuchsian groups and many other non-amenable groups.  The algebraic
eigenvalue properties derived in this paper can be summarized in the
following theorem.
\begin{theorem}[Corollary \ref{cor:notrans}, Theorems
   \ref{thm:free-maep}, \ref{thm:extension-maep}, \ref{fuchsian-maep}]
   Let $\sigma$ be an algebraic multiplier for the discrete group
   $G$, and let $A\in M_d(\Qbar(G,\sigma))$ be
   an operator acting on $l^2(G)^d$ by left multiplication
   twisted by $\sigma$, where $\Qbar$ denotes the field of algebraic
   numbers. Alternatively, consider $A$ to be a
   finite sum of magnetic translation operators
   $\sum_{g\in G} w_g R^\sigma_g$,
   where $w_g\in M_d(\Qbar)$.
   Then $A$ has only algebraic eigenvalues
   whenever
   \begin{enumerate}
   \item
      $G$ is finitely generated amenable, free or a surface group, or
   \item
      $\sigma$ is rational ($\sigma^n$=1 for some positive integer $n$)
      and $G\in\mathcal K$, where $\mathcal K$ is a class of
      groups containing free groups and amenable groups, and is closed
      under taking directed unions and extensions with amenable
      quotients.
   \end{enumerate}
\end{theorem}

The  case when $\sigma$ is a rational multiplier is established by
relating the spectrum of these
operators to the spectrum of untwisted operators on a finite covering
graph. The property follows from the class $\mathcal K$ having
the (untwisted) algebraic eigenvalue property (established in
section \ref{sec:algmult}, following \cite{dlmsy}), and from the
fact that $\mathcal K$ is closed under taking extensions with
cyclic kernel, as demonstrated in section \ref{sec:cyclic}.
We show in section \ref{sec:specproprational}
that these operators with rational weight function have no eigenvalues
that are Liouville transcendental whenever the group is
residually finite or more generally in a certain large class of groups
$\hat\calG$ containing $\mathcal K$. We also show that there is an
upper bound for the number of eigenvalues whenever the group satisfies
the Atiyah conjecture.

However, the case when $\sigma$ is an algebraic multiplier is
established in a significantly
different manner:  in addition to an approximation argument that
parallels that of \cite{dlmsy},
one also has to use new arguments that rely upon the geometry of closed
Riemann surfaces.

We also wish to highlight the
remarkable computation of Grigorchuk and \.Zuk \cite{GZ}, that is
recalled in Theorem \ref{GriZuk}. The computation explicitly lists the
dense set of eigenvalues of the Random Walk operator on the Cayley
graph of the lamplighter group --- all of these eigenvalues are algebraic
numbers, as predicted by results in this paper and in \cite{MY}, since
the lamplighter group is an amenable group.

Section 7 establishes an equality between the von Neumann spectral
density function of $A\in M_d(\C (G,\sigma))$ for arbitrary multiplier
$\sigma$, and the integrated density of states of $A$ with respect to
a generalized F\o{}lner exhaustion of $G$, whenever $G$ is a finitely
generated amenable group, or a surface group.

{\bf Acknowledgement} We would like to thank the referee for detailed
suggestions on improving the paper.

\section{Magnetic translations and the twisted group algebra}
\label{sec:magtrans}

The Harper operator is an example of an operator that can be
described as a sum of magnetic right translations. In this
section we will offer a brief description of the magnetic
translation operators, and observe that finite weighted sums
of these magnetic translations are unitarily equivalent to left
multiplication by matrices over a twisted group algebra.

Let $G$ be a discrete group and $\sigma$ be a {\em multiplier},
that is $\sigma\in Z^2(\Gamma;U(1))$
is a normalized $U(1)$-valued cocycle, which is a map from
$\Gamma \times \Gamma$ to $U(1)$ satisfying
\begin{gather}\label{cocycle}
   \sigma(b,c)\sigma(a,bc)=
   \sigma(ab,c)\sigma(a,b)\quad\forall a,b,c\in G
   \\
   \label{normalize}
   \sigma(1,g)=\sigma(g,1)= 1 \quad\forall g\in G
\end{gather}

Consider the Hilbert space $l^2(G)^d$ of square-integrable
$\C^d$-valued functions on $G$, where
$l^2(G)=\{h\colon G\to\C\,:\,\sum_{g\in G}\abs{h(g)}^2<\infty\}$.
The \emph{right magnetic translations} are then defined by
\[
   (R^\sigma_g f)(x)=f(xg)\sigma(x,g).
\]
Obviously, if $\sigma$ is a multiplier, so is $\overline{\sigma}$.
The right magnetic translations
commute with left magnetic translations $L^{\overline{\sigma}}_g$ as
follows from (\ref{cocycle}):
\begin{equation} \label{eqn:commuting_magnetic_translations}
   R^\sigma_g L^{\overline{\sigma}}_h
   =
   L^{\overline{\sigma}}_h R^\sigma_g\quad\forall g,h\in G
\end{equation}
where
\[
   (L^{\overline{\sigma}}_h f)(x)=f(h^{-1}x)\overline{\sigma(h,h^{-1}x)}.
\]

Consider now self-adjoint operators on $l^2(G)^d$ of the form
\begin{equation}
   \label{eqn:rt}
   A^\sigma = \sum_{g\in S} A(g) R^\sigma_g
\end{equation}
where $A(g)$ is a $d\times d$ complex matrix for each $g$, and
$S$ is a finite subset of $G$. The self-adjointness condition
is equivalent to demanding that the weights $A(g)$ satisfy
$A(g)^\ast=A(g^{-1}) \sigma(g,g^{-1})$.


These operators include as a special case the Harper operator and the  
DML
on the Cayley graph of $G$, where $S$ is the generating set
and $A(g)$ is identically $1$ for $g$ in $S$. For the Harper
operator of Sunada \cite{Sun} on a graph with finite fundamental
domain under the free action of the group $G$, one can construct
an operator of the form (\ref{eqn:rt}) which is unitarily
equivalent to the Harper operator, after identifying scalar
valued functions on the graph with $\C^n$-valued functions
on the group, where $n$ is the size of the fundamental domain.

Note that the operators of the form (\ref{eqn:rt}) are given explicitly
by the formula
$$
 (A^\sigma f)(x) = \sum_{g\in G} A(x^{-1}g)\sigma(x,x^{-1}g)f(g).
$$

For a given multiplier $\sigma$ taking values in
$\mathcal U(K)=K\cap U(1)$
for some subfield $K$ of the complex numbers,
one can also construct the twisted group algebra $K(G,\sigma)$
and examine $d\times d$ matrices $B^\sigma\in M_d(K(G,\sigma))$
acting on $l^2(G)^d$. Elements of $K(G,\sigma)$ are finite sums
$\sum a_g g$, $a_g\in K$ with multiplication given by
\[
   \left(\sum a_g g\right)\cdot \left(\sum b_g g\right)
   =
   \sum_{gh=k} a_g b_h \sigma(g,h) k.
\]
The action of $B^\sigma$ on a $\C^d$-valued function $f$ is
then given by this multiplication,
$$
   (B^\sigma f)(x) = \sum_{h\in G} B(xh^{-1})\sigma(xh^{-1},h)f(h)
$$
where $B(g)$ denotes the $d\times d$ matrix over $K$ whose elements
are the coefficients of $g$ in the elements of $B$. A straightforward
computation shows that 
\begin{equation}\label{eqn:lt}
   B^\sigma = \sum_{g\in S} B(g)L^\sigma_g,
\end{equation}
where $S$ is a finite subset of $G$.

The left and right twisted translations are unitarily equivalent
via the map $U^\sigma$,
\begin{gather*}
   (U^\sigma f)(x) = \sigma(x,x^{-1}) f(x^{-1}),
   \\
   U^\sigma L^\sigma_g = R^\sigma_g U^\sigma.
\end{gather*}
As such, operators of the form (\ref{eqn:rt}) and (\ref{eqn:lt}) 
will be unitarily equivalent if the coefficient matrices satisfy
$A(g)=B(g)$ for all $g\in G$.
Hereafter we will therefore concentrate on the latter picture,
noting that the results apply equally well to the case of operators
described as weighted sums of right magnetic translations.

By virtue of
(\ref{eqn:commuting_magnetic_translations}), the operators
$B^\sigma$ of the form (\ref{eqn:lt}) belong to the commutant
$B(l^2(G)^d)^{G,\sigma}$ of the set of magnetic translations
$\{R_g^{\overline{\sigma}}\,|\,g\in G\}$; the weak closure of
the set of operators of the form (\ref{eqn:lt}) is actually equal to 
this commutant as the theorem below shows.

The theorem itself is folklore, but we were not able to find the proof 
in the literature. 
In the special case of $G=\Z^2$, the details are spelt out in \cite{Sh}.
We will give a self-contained account, adapting the proof for the case
of trivial multiplier.

\begin{theorem}[Commutant theorem]
   \label{thm:commutant}
   The commutant of the right $\sigma$-translations on $l^2(G)$ is the
   von Neumann algebra generated by left $\bar\sigma$-translations on
   $l^2(G)$.
   
   Similarly, the commutant of the left $\sigma$-translations on
   $l^2(G)$ is the von Neumann algebra generated by right
   $\bar\sigma$-translations on $l^2(G)$.
\end{theorem}

\begin{proof}
   We present a proof of the second statement: the proof of the first
   statement is analogous.

   Let $W_{L,\sigma}$ be the von Neumann algebra generated by
   the set $S_{L,\sigma}=\{L_g^\sigma\,|\,g\in G\}$
   of left $\sigma$-translations, and
   $W_{R,\bar\sigma}$ be the von Neumann algebra generated by
   the set $S_{R,\bar\sigma}=\{R_g^{\bar\sigma}\,|\,g\in G\}$
   of right $\bar\sigma$-translations. We proceed
   by showing that $S_{L,\sigma}'=S_{R,\bar\sigma}''$
   (denoting the commutant by $'$) and then show that
   $S_{R,\bar\sigma}''=W_{R,\bar\sigma}$.

   An operator $C\in B(l^2(G))$ is determined by its components
   $C_{a,b}=(C\delta_b,\delta_a)=(C\delta_b)(a)$ for $a,b\in G$.
   Suppose $C\in S_{R,\bar\sigma}'$. In terms of components,
   one has that
   $(R_g^{\bar\sigma})_{a,b}
   =\delta_b(ag)\overline{\sigma(a,g)}
   =\delta_a(bg^{-1})\overline{\sigma(bg^{-1},g)}$,
   giving
   $
      (CR_g^{\bar\sigma})_{a,b}
      =C_{a,bg^{-1}}\overline{\sigma(bg^{-1},g)}$ and 
      $
      (R_g^{\bar\sigma}C)_{a,b}
      =\overline{\sigma(a,g)}C_{ag,b}$.
   $C$ commutes with $R_g^{\bar\sigma}$ for all $g$, and so
   substituting $bg$ for $b$ gives
   \[
   C\in S_{R,\bar\sigma}'
   \implies
   C_{a,b}=\overline{\sigma(a,g)}C_{ag,bg}\,\sigma(b,g)
   \quad\forall a,b,g\in G.
   \]

   Similarly for $D\in S_{L,\sigma}'$, noting that
   $(L_g^\sigma)_{a,b}
   =\delta_b(g^{-1}a)\sigma(g,g^{-1}a)
   =\delta_a(gb)\sigma(g,b)$,
   we have $(DL_g^\sigma)_{a,b}=D_{a,gb}\sigma(g,b)$ and
   $(L_g^\sigma D)_{a,b}=\sigma(g,g^{-1}a)D_{g^{-1}a,b}$,
   which after substituting $ga$ for $a$ gives
   \[
   D\in S_{L,\sigma}'
   \implies
   D_{a,b}=\overline{\sigma(g,a)}D_{ga,gb}\,\sigma(g,b)
   \quad\forall a,b,g\in G.
   \]

   Consider the product $CD$ for $C\in S_{R,\bar\sigma}'$
   and $D\in S_{L,\sigma}'$. In terms of components,
   \begin{align*}
   (CD)_{a,b}
   =
   \sum_{g\in G} C_{a,g}D_{g,b}
   &=
   \sum_{g\in G} \phi(a,g^{-1},b) C_{a(g^{-1}b),b}D_{a,(ag^{-1})b}
   \\
   &=
   \sum_{h\in G} \phi(a,a^{-1}hb^{-1},b) D_{a,h}C_{h,b},
   \end{align*}
   where
   $\phi(a,g^{-1},b)=
   \overline{\sigma(a,g^{-1}b)}\sigma(g,g^{-1}b)
   \overline{\sigma(ag^{-1},g)}\sigma(ag^{-1},b)$. However we
   can reduce the expression for $\phi$ by applying the cocycle
   identities to show that it is in fact identically equal to $1$:
   \begin{align*}
   \phi(a,k,b)
   &=
   \sigma(k^{-1},kb)\left(\sigma(ak,b)\overline{\sigma(a,kb)}\right)
   \overline{\sigma(ak,k^{-1})}
   \\
   &=
   \sigma(k^{-1},kb)\sigma(k,b)\overline{\sigma(a,k)\sigma(ak,k^{-1})}
   \\
   &=
   \sigma(k^{-1},k)\overline{\sigma(k,k^{-1})}
   \\
   &=1\qquad\forall a,k,b\in G.
   \end{align*}
   So $(CD)_{a,b}=(DC)_{a,b}$ for all $a,b\in G$,
   demonstrating that operators in $S_{L,\sigma}'$
   commute with those in $S_{R,\bar\sigma}'$.
   
   This gives the inclusion $S_{L,\sigma}'\subset S_{R,\bar\sigma}''$.
   The left $\sigma$-translations and
   right $\bar\sigma$-translations commute,
   so we also have that $S_{L,\sigma}\subset S_{R,\bar\sigma}'$
   and thus $S_{R,\bar\sigma}''\subset S_{L,\sigma}'$.
   Therefore $S_{L,\sigma}'=S_{R,\bar\sigma}''$.

   A calculation shows that the adjoint of $R^{\bar\sigma}_g$ is given by 
   $(R^{\bar\sigma}_g)^*=\sigma(g,g^{-1})R^{\bar\sigma}_{g^{-1}}$,
   and so operators that commute with the right $\bar\sigma$-translations
   must commute with their adjoints as well. So
   $S_{R,\bar\sigma}'=S_{R,\bar\sigma}^{*\prime}$,
   writing $S^\ast$ for the set of adjoints of elements of $S$.

   By the von Neumann double commutant theorem, the algebra generated
   by a set $S$ is given by $(S\cup S^\ast)''$.
   So
   $W_{R,\bar\sigma}=
   (S_{R,\bar\sigma}\cup S_{R,\bar\sigma}^\ast)''=
   S_{R,\bar\sigma}''=S_{L,\sigma}'$.
\end{proof}


We set the notation
$W^*_L(G,\sigma)=W_{L,\sigma}=S_{R,\bar\sigma}'$
for the left twisted group von Neumann algebra and
$W^*_R(G,\sigma)=W_{R,\sigma}=S_{L,\bar\sigma}'$
for the right twisted group von Neumann algebra.

The following is a corollary of the theorem.

\begin{corollary}
   \label{cor:commutant2}
   The commutant of the right $\sigma$-translations on $l^2(G)^d$ 
   is the von Neumann algebra $W_L^*(G, \bar\sigma) \otimes M_d(\C)$.

   Similarly, the commutant of the left $\sigma$-translations on $l^2(G)^d$ 
   is the von Neumann algebra $W_R^*(G, \bar\sigma) \otimes M_d(\C)$.
\end{corollary}

\begin{theorem}[Existence of trace]
   There is a canonical faithful, finite and normal trace on 
   the twisted group von Neumann algebras $W_L^*(G, \sigma)$ and
   $W_R^*(G, \bar\sigma)$ which is given by
   \begin{equation}
      \tr_{G,\sigma}(A)  = (A\delta_e,\delta_e).
   \end{equation}
   This trace is weakly continuous and can also be written as
   \begin{equation}
      \tr_{G,\sigma}(A) = (A\delta_g,\delta_g),\quad g\in G.
   \end{equation}
\end{theorem}

\begin{proof}
   It is clear that $\tr_{G,\sigma}$ is linear, finite and weakly
   continuous (hence normal).  Now if $A\in W^*_L(G, \sigma)$, then
   \begin{equation}
      (A\delta_g,\delta_g) =  A_{g, g} =
      \sigma(g, h) A_{gh, gh} \overline{\sigma(g, h)} = A_{gh, gh}
      =(A\delta_{gh},\delta_{gh}) 
   \end{equation}
   for all $h\in G$. In particular, every diagonal entry of the matrix
   of $A$ is equal to $\tr_{G,\sigma}(A)$.
   \\
   
   If $A$ is a self-adjoint operator in $W^*_L(G, \sigma)$ such that
   $\tr_{G,\sigma}(A)=0$, then $(A\delta_g,\delta_g)=0$ for all
   $g\in G$.  But then due to the Cauchy-Schwarz inequality
   $|(Af_1,f_2)|^2\le (Af_1,f_1)(Af_2,f_2),$ $\ f_1,f_2\in l^2(G)$,
   we deduce that $(A\delta_g,\delta_h)=0$ for all $g,h\in G$, which
   implies that $A=0$. Therefore
   $\tr_{G,\sigma}$ is faithful.
   \\

   It remains to prove that $\tr_{G,\sigma}$ is a trace.  That is, 
   \begin{equation}
      \tr_{G,\sigma}(AB) = \tr_{G,\sigma}(BA),\quad A,B \in W^*_L(G, \sigma)
   \end{equation}
   Since $\tr_{G,\sigma}$ is linear and weakly continuous it is
   sufficient to consider the case when $A=L^\sigma_g$ and
   $B=L^\sigma_h$ for all $ g,h\in G$.  We compute,
   \begin{align*}
      \tr_{G,\sigma}(L^\sigma_g L^\sigma_h)
      & = (L^\sigma_g L^\sigma_h\delta_e,\delta_e)
      = \sigma(g, h) (L^\sigma_{gh}\delta_e,\delta_e)
      = \sigma(g, h) (\delta_{gh},\delta_e)
      \\[+7pt]
      & =
      \begin{cases}
         \sigma(g, h) & \text{if } gh=e, \\
         0            & \text{otherwise}.
      \end{cases}
   \end{align*}
   Similarly, 
   $$
   \tr_{G,\sigma}(L^\sigma_h L^\sigma_g)=
   \begin{cases}
      \sigma(h, g) & \text{if } hg=e, \\
      0            & \text{otherwise}.
   \end{cases}
   $$
   By the cocycle identity \eqref{cocycle} with $a=h^{-1}$, $b=h$
   and $c=h^{-1}$, we see that
   $$
   \sigma(h, h^{-1}) = \sigma(h^{-1}, h) \qquad \forall h\in G.
   $$
   Therefore
   $ \tr_{G,\sigma}(L^\sigma_g L^\sigma_h) =
   \tr_{G,\sigma}(L^\sigma_h L^\sigma_g)$
   for all $g,h\in G$ as required.
   The argument for $W_R^*(G, \bar\sigma)$ is identical.
\end{proof}

Now the matrix algebra $M_d (\C)$ has the canonical trace $\Tr$ given
by the sum of the diagonal coefficients of a matrix. Then the tensor
product $\tr_{G,\sigma} \otimes \Tr$ is a trace on
$W_L^*(G, \sigma)\otimes M_d(\C)$ and on
$W_R^*(G, \bar\sigma) \otimes M_d( \C)$ which we denote by
 $\tr_{G,\sigma}$ for simplicity.

\begin{corollary}\label{cor:trace}
   There is a canonical faithful, finite and normal trace on the
   twisted group von Neumann algebras $W_L^*(G, \sigma) \otimes M_d(
   \C)$ and on $W_R^*(G, \bar\sigma) \otimes M_d( \C)$ which is given
   by
   \begin{equation}
      \tr_{G,\sigma}(A)  = \sum_{j=1}^d (A_{j,j} \delta_e,\delta_e).
   \end{equation}
   This trace is weakly continuous and can also be written as
   \begin{equation}
      \tr_{G,\sigma}(A) =
      \sum_{j=1}^d (A_{j,j} \delta_g,\delta_g),\quad g\in G.
   \end{equation}
\end{corollary}

In the corollary above, we interpret the elements of the tensor
product $W_L^*(G, \sigma) \otimes M_d(\C) = M_d(W_L^*(G,\sigma))$
as $d\times d$ matrices with entries in $W_L^*(G, \sigma)$ etc.

Suppose that $\calA$ is a von Neumann algebra of algebras
of operators acting on a Hilbert space $\calH$. A subspace $U$
of $\calH$ is termed \emph{affiliated} if the corresponding
orthogonal projection $P_U$ onto the closure of $U$ belongs to $\calA$.
A necessary and sufficient condition for affiliation
is that the subspace be invariant under the action of operators
in the commutant $\calA'$ of $\calA$. We will write $U\eta \calA$ to
indicate that the subspace $U$ is affiliated to the algebra $\calA$.

Given a trace $\tau$ on $\calA$, the von Neumann dimension $\dim_\tau$
of an affiliated subspace is defined to be the trace of $P_U$.
We will use the following properties of the von Neumann
dimension (see for example \cite{Sh2}, Section 2.6 Lemma 2,
Section 2.26).
\begin{lemma}
  \label{l:vndim1}
  Let $\calH$ be a Hilbert space, $\calA$ a von Neumann
  algebra of operators on $\calH$ with (normal, faithful
  and semi-finite) trace $\tau$ 
  and von Neumann dimension $\dim_\tau$, and let
  $L,N\subset \calH$ be affiliated subspaces. Then,
  \begin{enumerate}
    \item
      $\dim_\tau L=0$ implies $L=\{0\}$,
    \item
      $L\subseteq N$ implies $\dim_\tau L\leq\dim_\tau N$,
    \item
      if $A\in\calA$ is an almost
      isomorphism of $L$ and $N$, that is, $ker A\cap L=\{0\}$
      and the set $A(L)$ is dense in $N$, then
      $\dim_\tau L=\dim_\tau N$.
    \end{enumerate}
\end{lemma}
The following is an immediate consequence.
\begin{lemma}
  \label{l:vndim2}
  Let $\calH$, $\calA$, $\tau$ be as in Lemma \ref{l:vndim1}.
  If $L$ is an affiliated subspace of $\calH$ with corresponding
  projection $P_L$, then
  \[
  \dim_\tau L=\dim_\tau \ker A|_L+\dim_\tau \overline{\im A|_L}
  =\dim_\tau (\ker A\cap L) + \dim_\tau \overline{\im AP_L}.
  \]
\end{lemma}
\begin{proof}
  This follows by noting that $A$ gives an almost isomorphism
  from the orthogonal complement of its kernel in $L$ to
  the closure of its image on $L$.
\end{proof}

Hereafter we will use $\dim_{G,\sigma}$ to refer to the von Neumann
dimension associated with the trace $\tr_{G,\sigma}$ on the algebra
$W^*_L(G, \sigma)\otimes M_d( \C)$. In the case that $\sigma$ is trivial, this
algebra becomes the von Neumann algebra of $G$-equivariant operators
$B(l^2(G)^d)^G$, and we refer to the trace and dimension
by $\tr_G$ and $\dim_G$ respectively. To make the notation more
suggestive we will also write
$W^*_L(G, \sigma)\otimes M_d( \C)=B(l^2(G)^d)^{G,\sigma}$.

Two multipliers $\sigma$ and $\sigma'$ in
$\Z^2(G,\mathcal U(K))$ are \textit{cohomologous},
written $\sigma\sim\sigma'$, if they belong to the same cohomology
class in $H^2(G,\mathcal U(K))$. It follows that $\sigma\sim\sigma'$
if and only if there exists a map $s:G\to \mathcal U(K)$ such that
\begin{equation}
   \label{eq:cohomologous}
   \sigma(g,f)=s(g)s(h)\overline{s(gh)}\sigma'(g,h)
   \quad\forall g,h\in G.
\end{equation}
The map $s$ gives rise to a unitary equivalence between operators
in $M_d(K(G,\sigma))$ and $M_d(K(G,\sigma'))$.

\begin{lemma}
   \label{l:cohom1}
   Let $\sigma$ and $\sigma'$ be cohomologous multipliers in
   $\Z^2(G,\mathcal U(K))$.
   Then for every $A$ in $M_d(K(G,\sigma))$ acting on $l^2(G)^d$
   there is a canonically determined $A'$ in $M_d(K(G,\sigma'))$ such that
   $A$ and $A'$ are unitarily equivalent.
\end{lemma}

\begin{proof}
   Let $s:G\to \mathcal U(K)$ be the map as in (\ref{eq:cohomologous}),
   such that $\sigma(g,h)=s(g)s(h)\overline{s(gh)}\sigma'(g,h)$ 
   for all $x,y\in G$. Writing $A(g)\in M_d(K)$ for the matrix
   of coefficients of $g$ in $A$, as in (\ref{eqn:lt}), one has
   \begin{align*}
     (Af)(g) &=\sum_{h\in G} A(gh^{-1})f(h)\sigma(gh^{-1},h)
     \\
     &=\sum_{h\in G}
     A(gh^{-1})f(h)s(gh^{-1})s(h)s(g)^{-1}\sigma'(gh^{-1},h).
   \end{align*}
   Let $S$ be the unitary operator on $l^2(G)^d$ given by multiplication
   by $s$: $Sf(g)=s(g)f(g)$. Then letting $A'(g)=s(g)A(g)$, one has
   \begin{align*}
     (SAf)(g)
     &=\sum_{h\in G} A'(gh^{-1})f(h)s(h)\sigma'(gh^{-1},h).
     \\
     &=(A'Sf)(g)
   \end{align*}
   for all $g\in G$, $f\in l^2(G)^d$. That is, $A$ and $A'$ are   
   unitarily equivalent.
\end{proof}

It is sometimes convenient to consider only the case when the multiplier
satisfies $\sigma(g,g^{-1})=1$ for all $g$ in $G$. The following lemma
shows that there is such a multiplier in every cohomology class when
the subfield $K$ is algebraically closed.

\begin{lemma}
   \label{l:nicemultiplier}
   Suppose $K$ is an algebraically closed subfield of $\C$.
   Then any multiplier $\sigma\in Z^2(G,\mathcal U(K))$
   is cohomologous to a multiplier $\sigma'$ such that
   $\sigma'(g,g^{-1})=1$ for all $g\in G$.
\end{lemma}
\begin{proof}
   By the cocycle identity, $\sigma(g,g^{-1})=\sigma(g^{-1},g)$ for
   all $g\in G$. Choose $s:G\to\mathcal U(K)$ such that
   $s(g)=s(g^{-1})$ and $s(g)^2=\sigma(g,g^{-1})$, for example
   by setting $s(g)=e^{i\theta/2}$ when $\sigma(g,g^{-1})=e^{i\theta}$,
   for $\theta\in [0,2\pi).$
   The image of $s$ lies in $\mathcal U(K)$ due to $K$ being
   algebraically closed.

   Let $\sigma'$ be the cohomologous multiplier given by $s$,
   according to the formula (\ref{eq:cohomologous}). Then
   \begin{align*}
      \sigma'(g,g^{-1})
      &=\overline{s(g)s(g^{-1})}s(1)\sigma(g,g^{-1})
      \\
      &=\overline{\sigma(g,g^{-1})}\sigma(g,g^{-1})
      \\
      &=1,\qquad\forall g\in G.
   \end{align*}
\end{proof}

\section{Algebraic eigenvalue property}
\label{sec:aep}

The algebraic eigenvalue property for groups
was introduced in \cite{dlmsy}.
We recall the definition here, and present a class of groups
$\calK$ for which the algebraic eigenvalue property holds.
We then define a similar property describing the
eigenvalues for matrix operators over
the twisted group ring, as described in section \ref{sec:magtrans},
termed the $\sigma$-multiplier algebraic eigenvalue property.

Thus recall the following definition from \cite{dlmsy}, where
$\Qbar$ denotes the set of complex algebraic numbers.

\begin{definition}[4.1 of \cite{dlmsy}]\label{aep}
   A discrete group $G$ has the
   \emph{algebraic eigenvalue property}, if for every $d\times d$
   matrix $A\in M_d(\Qbar G)$ the eigenvalues of $A$, acting
   on $l^2(G)^d$, are algebraic numbers.
\end{definition}

Note that operators without point spectrum satisfy the criterion
in the vacuous sense.

The trivial group has the
algebraic eigenvalue property, since the eigenvalues are the zeros of
the characteristic polynomial. The same is true for every finite
group. More generally, if $G$ contains a subgroup $H$ of finite index,
and $H$ has the algebraic eigenvalue property, then the same is true
for $G$. And if $G$ has the algebraic eigenvalue property and $H$ is a
subgroup of $G$, then $H$ also has the algebraic eigenvalue property.

In section 4 of \cite{dlmsy} it was shown that the algebraic
eigenvalue property holds for all amenable groups and for all groups
in Linnell's class $\calC$, which is the smallest class of groups
containing all free groups and which is closed under extensions
with elementary amenable quotient and under directed unions. This
motivates the definition of the class $\calK$, a larger class
which contains these groups, for which the algebraic eigenvalue
property can be shown to hold.

\begin{definition}
   \label{defn:calk}
   The class $\calK$ is the
   smallest class of groups containing free groups and
   amenable groups, which is closed under taking extensions with   
   amenable quotient, and under taking directed unions.
\end{definition}

\begin{remark}
   It is clear that the class $\calK$ contains
   every discrete amenable group and every group in Linnell's class
   $\calC$.

   Recall that the class of elementary amenable groups is the smallest
   class of groups containing all cyclic and all finite groups and which
   is closed under taking group extensions and directed unions.
   As such then $\calK$ is a strictly larger class of groups than
   $\calC$,
   as it contains amenable groups which are not elementary amenable,
   such as the example presented by Grigorchuk in \cite{G83}.
\end{remark}

\begin{remark}\label{surface}
   Every subgroup of infinite index in a surface group $\Gamma$ is a
   free group.  Here $\Gamma$ is the fundamental group of a compact
   Riemann surface of genus $g>1$.  This follows from the fact that such
   groups are fundamental groups of an infinite cover of the base
   surface and from the general fact that the fundamental group of a
   noncompact surface is free (see \cite{AS} Chapter 1, \S~7.44 and \S~8.)
   Since we have the exact sequence
   \[
     1\to F\to \Gamma \to \mathbb Z^{2g}\to 1
   \]
   where $F$ is a free group by the argument above and the free
   abelian group $ \mathbb Z^{2g}$ is an elementary amenable group, we
   deduce that the surface group $\Gamma$ belongs to the class $\calC$,
   and hence also to the class $\calK$.
\end{remark}

\begin{remark}\label{fuchsian}
   Let $\Gamma$ be a cocompact Fuchsian group, namely $\Gamma$ is a
   discrete subgroup of ${\bf SL}(2, \mathbb R)$ such that the quotient
   space $\Gamma\backslash {\bf SL}(2, \mathbb R)$ is compact.
   Then there is a torsion-free subgroup $G$ of $\Gamma$ of
   finite index such that $G$ is the fundamental group of a compact
   Riemann surface of genus greater than one. By Remark \ref{surface} above,
   $G$ is in $\calK$, and since $\Gamma/G$ is a finite group, it is
   amenable. Therefore $\Gamma$ is also in $\calK$.
\end{remark}

\begin{remark}\label{modular}
  Consider the the modular group ${\bf SL}(2, \mathbb Z) $.  Then it
  is well known that there is a congruence subgroup $\Gamma(N)$ of
  finite index in ${\bf SL}(2, \mathbb Z)$ that is isomorphic to a
  free group.  We conclude by the arguments in Remarks \ref{surface}
  and \ref{fuchsian}, that the modular group and all of the congruence
  subgroups are in $\calK$.
\end{remark}

\begin{theorem}
   \label{thm:Kalg}
   Every group in $\calK$ has the algebraic eigenvalue property.
\end{theorem}

The proof of this theorem closely follows the argument
in \cite{dlmsy} for $\calC$, and we leave the details to
section \ref{sec:classK}.

\begin{remark}\label{letf-right}
   Results of \cite{dlmsy} were formulated for operators of the form
   \begin{equation}\label{left}
     B = \sum_{g\in S} B(g) L_g
   \end{equation}
   acting on $l^2(G)^n$ where $S\subset G$ is a finite subset,
   $A(g)$ is an $n\times n$ complex matrix and $L_g$ denotes the
   \emph{untwisted} left translation on $l^2(G)$. Since
   $x\mapsto x^{-1}$ induces a unitary transformation on $l^2(G)$
   that conjugates the right translation $R_g$ with $L_g$,
   we see that (\ref{left}) is unitarily equivalent to
   \begin{equation}\label{right}
     \sum_{h\in S} B(h) R_h.
   \end{equation}
   It follows that all results of \cite{dlmsy} concerning spectral
   properties of operators (\ref{left}) apply to operators of the form
   (\ref{right}) equally well.
\end{remark}
Suppose now we have an operator $A\in M_d(\Qbar(G,\sigma))$
acting on $l^2(G)^d$ by left twisted multiplication, as described
in section \ref{sec:magtrans}, where $\Qbar(G,\sigma)$ is the
twisted group algebra over the algebraic numbers $\Qbar$ with
multiplier $\sigma$. For a fixed $\sigma$, one can ask
if any such $A$ can have transcendental eigenvalues.

\begin{definition}
   A discrete group $G$ is said to have the
   \emph{$\sigma$-multiplier algebraic eigenvalue property},
   if for every matrix $A\in M_d(\Qbar(G,\sigma))$,
   regarded as an operator on $l^2(G)^d$, the eigenvalues of $A$
   are algebraic numbers, where
   $\sigma \in Z^2(G, \mathcal U(\Qbar))$ is
   an algebraic multiplier, and $\mathcal U(\Qbar)$
   denotes the unitary elements of the field of algebraic numbers.
\end{definition}

An immediate consequence of Lemma \ref{l:cohom1} is that for
a given group $G$, the $\sigma$-multiplier algebraic eigenvalue
property depends only on the cohomology class of $\sigma$.

\begin{corollary}
   \label{cor:cohom-maep}
   Suppose $G$ has the $\sigma$-multiplier algebraic eigenvalue
   property. Then $G$ has the $\sigma'$-multiplier algebraic eigenvalue
   property for any $\sigma'\sim\sigma$ in $Z^2(G,\mathcal U(\Qbar))$.
\end{corollary}
\begin{proof}
   Any $A'\in M_d(\Qbar(G,\sigma'))$ is unitarily equivalent to some
   $A\in M_d(\Qbar(G,\sigma))$ by Lemma \ref{l:cohom1}, and so has
   only algebraic eigenvalues.
\end{proof}

In the following sections \ref{sec:specproprational} and
\ref{sec:findistinct} we investigate the situation when
$\sigma$ is rational, that is, when $\sigma^n=1$ for some $n$.
In particular it is shown that every group in $\calK$ has
the $\sigma$-multiplier algebraic eigenvalue property
when $\sigma$ is rational.

The case of more general algebraic multipliers is
discussed in section \ref{sec:algmult}.

\section{Spectral properties with rational $\sigma$}
\label{sec:specproprational}

Suppose the weight function $\sigma$ is rational with $\sigma^r=1$,
and let $G^\sigma$ be the extension of $G$ by $\Z_r$
as follows,
\begin{equation}
   \label{e:Gammasigma}
   \begin{gathered}
     1\longrightarrow\Z_r\longrightarrow G^\sigma
     \longrightarrow G\longrightarrow 1
     \\
     (z_1,g_1)\cdot(z_2,g_2)=
     (z_1 z_2 \sigma(g_1,g_2),g_1 g_2)
   \end{gathered}
\end{equation}
regarding $\Z_r$ as a (multiplicative) subgroup of $U(1)$.
One can then relate the spectrum of an operator
$A^\sigma\in M_d(K(G,\sigma))$ acting on $l^2(G)^d$ as in \eqref{eqn:lt}
to that of an associated operator $\tilde A$ on $l^2(G^\sigma)^d$.

Define a map $\Psi:\,M_d(K(G,\sigma))\to M_d(KG^\sigma)$ as follows.
For $A^\sigma\in M_d(K(G,\sigma))$ with matrices of coefficients
$A(g)\in M_d(K)$, let $\tilde A=\Psi(A^\sigma)$ be given by
\begin{equation}
   \label{e:Atilde}
   \tilde{A}(z,g)=
   \begin{cases}
     A(g) &\text{ if $z=1$,} \\
     0    &\text{ otherwise,}
   \end{cases}
\end{equation}
acting on $l^2(G^\sigma)^d$ by left multiplication.

Consider the map $\xi:l^2(G)^d\to l^2(G^\sigma)^d$ given by
$(\xi f)(z,g) = \overline{z}f(g)$. Then $\frac{1}{\sqrt r}\xi$ is an
isometry from $l^2(G)^d$ to the closed subspace $R$ of
$l^2(G^\sigma)^d$ where $R=\{f|\,f(z,g)=\overline{z}f(1,g)\;
\forall (z,g)\in G^\sigma\}$. By \eqref{e:Gammasigma},
$(1,g)^{-1}\cdot (z,h)=(z\,\overline{\sigma(g,g^{-1}h)},g^{-1}h)$ and so
\begin{align*}
   (\tilde A\xi f)(z,h)
   &=
   \sum_{(z',g)\in G^\sigma}
   \tilde{A}(z',g)(\xi f)((z',g)^{-1}\cdot(z,h))
   \\
   &=
   \sum_{g\in G} A(g)(\xi f)(z\overline{\sigma(g,g^{-1}h)},g^{-1}h)
   \\
   &=
   \sum_{g\in G} A(g) f(g^{-1}h) \sigma(g,g^{-1}h) \overline{z}
   \\
   &=
   (\xi A^\sigma f)(z,h),
\end{align*}
for all $(z,h)\in G^\sigma$, and thus
\begin{equation}
   \label{e:upreln}
   \Psi(A^\sigma)\xi=\xi A^\sigma\quad\forall A^\sigma\in M_d(K(G,\sigma)).
\end{equation}
$A^\sigma$ is therefore unitarily equivalent to the restriction to the
subspace $R$ of the operator $\tilde A$. 

\begin{lemma}
  Let $A$ be a bounded linear operator on a Hilbert space $H$,
  such that $\im A\vert_R\subset R$ for a closed subspace $R$ of $H$. Then
  regarding $A\vert_R$ as an operator on $R$,
  $\specpoint A\vert_R\subset\specpoint A$ and $\spec A\vert_R\subset\spec A$.
\end{lemma}
\begin{proof}
  Any eigenvector in $R$ is an eigenvector in $H$ and so the inclusion of
  point spectrum is immediate.

  Suppose $\lambda\not\in\spec A$. Then $\lambda\not\in\specpoint A\vert_R$
  and $\im (A-\lambda)\vert_R$ is dense in $R$. Let $B$ be the inverse
  of $A-\lambda$. For $u\in R$ one can find a convergent net
  $u_\alpha\to u$ with $u_\alpha=(A-\lambda)u_\alpha'$ for $u_\alpha'$ in $R$.
  Applying $B$ gives $u_\alpha'\to Bu$, but $R$ is closed, and so $Bu$ is
  in $R$ and $u$ is in the image of $(A-\lambda)$.
  Therefore $(A-\lambda)\vert_R$ has inverse $B\vert_R$ and
  $\lambda\not\in\spec A\vert_R$.
\end{proof}

We therefore have spectral inclusions for $A^\sigma$ and $\tilde A$.

\begin{proposition}
   \label{prop:specinclcayley}
   \label{prop:ptspecinclcayley}
   Let $A^\sigma\in M_d(K(G,\sigma))$ be a bounded linear operator
   on $l^2(G)^d$ as in \eqref{eqn:lt}, and suppose $\sigma$ is rational.
   Let
   $\tilde A=\Psi(A^\sigma):\,l^2(G^\sigma)^d \to l^2(G^\sigma)^d$
   be the associated $G^\sigma$-equivariant operator as described above.
   Then
   \begin{align}
     \label{e:specinclcayley}
     \spec A^\sigma&\subseteq \spec \tilde A,
     \intertext{and}
     \label{e:ptspecinclcayley}
     \specpoint A^\sigma&\subseteq \specpoint \tilde A.
   \end{align}
\end{proposition}

The following is an easy corollary.

\begin{corollary}
   Let $\tilde A$ and $A^\sigma$ be as described in Proposition
   \ref{prop:specinclcayley}.  Then any interval $(a, b)$ that is a gap
   in the spectrum of $\tilde A$ is also contained in a gap of the
   spectrum of $A^\sigma$.
\end{corollary}

In section \ref{sec:cyclic} we prove the following result concerning
the class of groups $\calK$ introduced in section \ref{sec:aep}.

\begin{proposition}\label{prop:cyclic}
   The class $\calK$ of groups is closed under taking extensions
   with cyclic kernel.
\end{proposition}

Therefore, by Theorem \ref{thm:Kalg}, the groups in this class all
have the $\sigma$-multiplier algebraic eigenvalue property for
rational $\sigma$.

\begin{corollary}[Absence of eigenvalues that are transcendental
   numbers]
   \label{cor:notrans}
   Any $A^{\sigma} \in M_d(\Qbar(G,\sigma))$ has only
   eigenvalues that are algebraic numbers, whenever $G\in \calK$ and
   $\sigma$ is a rational multiplier on $G$.
\end{corollary}

\begin{proof}
   Let $G^\sigma$ be the central extension of the group
   $G$ in the class $\calK$, where $G^\sigma$ is
   defined in (\ref{e:Gammasigma}), and
   let $\tilde A$ be the operator on $l^2(G^\sigma)^d$ associated with
   $A^\sigma$ as defined in (\ref{e:Atilde}).
   By Proposition \ref{prop:cyclic}, any central extension of
   $G$ by a cyclic group $\Z_r$ is also in the class
   $\calK$, therefore the group $G^\sigma$ where
   $\sigma$  is in $\calK$.
   By Theorem \ref{thm:Kalg}, we know that every group
   $G$ in
   the class $\calK$ has the algebraic eigenvalue property
   and so $\tilde A$
   has only algebraic eigenvalues.
   $A^\sigma$ therefore has only algebraic eigenvalues by Proposition
   \ref{prop:ptspecinclcayley}.
\end{proof}

Recall the  definition of the following class of groups from
\cite{dlmsy}.

\begin{definition}
  Let $\calG$ be the smallest class of groups which contains the   
trivial
   group and is closed under the following processes:
   \begin{enumerate}
   \item If $H\in\calG$ and $G$ is a generalized amenable
     extension of $H$, then $G\in\calG$.
   \item If $H\in \calG$ and $U\subgroup H$, then $U\in\calG$.
   \item If $G=\lim_{i\in I}G_i$ is the direct or inverse limit of a
     directed system of groups $G_i\in \calG$,
     then  $G\in\calG$.
   \end{enumerate}
\end{definition}

We have the inclusion $\calC\subset\calK\subset\calG$,
and in particular the class $\calG$ contains all amenable groups,
free groups, residually finite groups, and residually amenable
groups. Consider the subclass of groups $\hat \calG$ defined as
$$
\hat\calG = \left\{G\in \calG: \hat G \in \calG \;\; \forall
   \text{ $\Z_r$-extensions $\hat G$ of $G$}
\right\}.
$$
By the results of section \ref{sec:cyclic}, we have the inclusion
$\calC\subset\calK\subset\hat\calG\subset\calG$.

\begin{corollary}[Absence of eigenvalues that are Liouville
  transcendental numbers]
  Any self-adjoint $A^\sigma\in M_d(\Qbar(G,\sigma))$ does not have
  any eigenvalues that are Liouville transcendental numbers, whenever
  $G\in\hat\calG$ and
  $\sigma$ is a rational multiplier on $G$.
\end{corollary}

\begin{proof}
   Any operator $A^\sigma\in M_d(\Qbar(G,\sigma))$
   is self adjoint if and only if $A(g)^*=A(g^{-1})\sigma(g,g^{-1})$
   for all $g\in G$. By Lemmas \ref{l:cohom1} and \ref{l:nicemultiplier},
   there exists an $A'\in M_d(\Qbar(G,\sigma'))$
   such that $A^\sigma$ and $A'$ are unitarily equivalent
   and where $\sigma'(g,g^{-1})=1$ for all $g\in G$.
   With $A^\sigma$ being self adjoint, we have that $A'$
   is self adjoint and thus that $A'(g)^*=A'(g^{-1})$ for all
   $g\in G$.
   By the construction of Lemma \ref{l:nicemultiplier}, if
   $\sigma$ is a rational multiplier with $\sigma^r=1$, then
   $\sigma'$ is also rational, with ${\sigma'}^{2r}=1$.

   Let $\tilde A=\Psi(A')\in M_d(\Qbar(G^{\sigma'}))$,
   as in (\ref{e:Atilde}), where $G^\sigma$ is the central
   extension of $G$ as described in (\ref{e:Gammasigma}).
   In terms of matrices of coefficients, $\tilde A(z,g)=\delta_1(z)A(g)$
   for $(z,g)\in G^{\sigma'}$. As $\sigma'(g,g^{-1})=1$,
   \[
   \tilde A((z,g)^{-1}) = \tilde A(z^{-1}\sigma(g,g^{-1}),g^{-1})
   = \tilde A(z^{-1},g^{-1})
   \]
   and
   \[
   \tilde A(z^{-1},g^{-1})
   = \delta_1(z^{-1})A'(g^{-1})
   = \delta_1(z)A'(g)^\ast
   = \tilde A(z,g)^*,
   \]
   showing that $\tilde A$ is self-adjoint.

   $\hat\calG$ is closed under extensions by finite cyclic groups,
   and so the group $G^{\sigma'}$ is in the  class
   $\calG$. Applying Theorem 4.15 of \cite{dlmsy}, it follows that
   $\tilde A$ does not have any eigenvalues that are Liouville
   transcendental numbers.
   Then by Proposition \ref{prop:ptspecinclcayley},
   $A'$ and thus $A^\sigma$ do not have
   any eigenvalues that are Liouville transcendental numbers.
\end{proof}

\section{On the finiteness of the number of distinct eigenvalues}
\label{sec:findistinct}

We deal here with the following situation: $G$ is a discrete group and
$A\in M_d(\Qbar G)$.  Then $A$
induces a bounded linear operator $A\colon l^2(G)^d\to l^2(G)^d$ by
left convolution (using the canonical left $G$-action on $l^2(G)$),
which commutes with the right $G$-action.

Let ${\pr}_{\ker A}\colon l^2(G)^d\to l^2(G)^d$ denote the orthogonal
projection onto $\ker A$. Recall that the von Neumann dimension of
$\ker A$ is defined as
\[
   \dim_G(\ker A): =\tr_G(\pr_{\ker A}):=
   \sum_{i=1}^d \langle{\pr_{\ker A} e_i,e_i}\rangle_{l^2(G)^d},
\]
where $e_i\in l^2(G)^d$ is the vector with the trivial element of
$G\subset l^2(G)$ at the $i^\th$-position and zeros elsewhere.

Let $G$ be a discrete group.
Let ${\fin}(G)$ denote the additive subgroup of
$\mathbb Q$ generated by the inverses of
the orders of the finite subgroups of $G$. Note that
$\fin(G)=\mathbb Z$ if and only if $G$ is torsion free and
$\fin(G)$ is discrete in $\mathbb R$ if and only if
orders of finite subgroups of $G$ are bounded.
Recall the following definition.

\begin{definition}\label{def:AtiyahCon}
   A discrete group $G$ is said to fulfill the
   \emph{strong Atiyah conjecture}
   if the orders of the finite subgroups of $G$ are bounded
   and
   \[
     \dim_G (\ker A)\in\fin(G) \qquad\forall A\in
     M_d(\Qbar G);
   \]
   where $\ker A$ is the kernel of the induced map $A\colon l^2(G)^d\to
   l^2(G)^d$.
\end{definition}

Linnell proved the strong Atiyah conjecture if $G$ is such that the
orders of the finite subgroups of $G$ are bounded and $G\in\calC$,
where $\calC$ is Linnell's class of groups that is defined just below
Definition \ref{aep}.

\begin{theorem}[\cite{Lin}]\label{theo:algebraicstrongAt}
   If $G\in\calC$ is such that the orders of the finite subgroups of $G$
   are bounded, then the strong Atiyah conjecture is true.
\end{theorem}

In \cite{dlmsy}, Linnell's results were generalized to a larger class
$\calD$ of groups, but these groups are all torsion-free and therefore
our results do not apply to them.

The following is an easy corollary of the relationship between
$A^\sigma$ and $\tilde A$, cf.\ (\ref{e:Gammasigma}).

\begin{corollary}[Finite number of distinct eigenvalues]
   \label{finiteev}
   Any self-adjoint $A^{\sigma} \in M_d(\Qbar(G,\sigma))$ has only a
   finite number of distinct eigenvalues whenever $G \in \calC$ is such
   that the orders of the finite subgroups of $G$ are bounded, and
   $\sigma$ is a rational multiplier on $G$.
\end{corollary}

\begin{proof}
   Let $G^\sigma$ be the central extension of the group $G$ in
   $\calC$, where $G^\sigma$ is defined in (\ref{e:Gammasigma}), and
   let $\tilde A$ be the operator on $l^2(G^\sigma)$ associated with
   $A^\sigma$ as defined in (\ref{e:Atilde}).  By Remark
   \ref{linnell}, the group $G^\sigma$ as defined in
   (\ref{e:Gammasigma}) is also in $\calC$. Clearly, the orders of
   finite subgroups of $G^\sigma$ are bounded as well. Thus the
   Theorem \ref{theo:algebraicstrongAt} above applies to $G^\sigma$,
   and so the dimensions of eigenspaces of $\tilde A$ are in the
   discrete, closed subgroup $\fin(G)$ of $\R$ and are therefore
   bounded away from zero. It follows that $\tilde A$ can have at most
   finitely many eigenvalues.  The conclusion now follows from
   (\ref{e:ptspecinclcayley}).
\end{proof}

\begin{remark}
   Our results do not just apply to operators acting on scalar valued
   functions but also to vector valued functions. In this case there
   are many examples where eigenvalues exist. For instance, for the
   combinatorial Laplacian on $L^2$, degree zero cochains of a
   covering space, zero is never an eigenvalue, whereas it is very
   common for it to be an eigenvalue on $L^2$ cochains of positive
   degree. This is the case whenever the Euler characteristic of the
   base is nonzero, which follows for instance from Atiyah's $L^2$
   index theorem for covering spaces \cite{At}, and Dodziuk's theorem
   on the combinatorial invariance of the $L^2$ Betti numbers,
   \cite{Do}.
\end{remark}


Let $H$ denote the lamplighter group, namely $H$ is the wreath product
of $\mathbb Z_2$ and $\mathbb Z$. Then there is the following
remarkable computation of Grigorchuk and \.Zuk, Theorem 2 and Corollary
3, \cite{GZ}.

\begin{theorem}\label{GriZuk}
   Let $A:= t+at+t^{-1}+(at)^{-1}\in \mathbb Z H$ be a multiple of the
   Random Walk operator of $H$. Then $A$, considered as an operator on
   $l^2(H)$, has eigenvalues
   \begin{equation}
     \left\{ 4\cos\left(\frac{p}{q}\pi\right)\mid
       p\in\mathbb Z , q=2,3,\dots\right\}.
   \end{equation}
   The $L^2$-dimension of the corresponding eigenspaces is
   \begin{equation}
     \dim_H \ker\left(A-4\cos\left(\frac{p}{q}\pi\right)\right) =
     \frac{1}{2^q-1}\qquad \text{if }p,q\in\mathbb {Z},
     \ q\ge 2, \text{ with }(p,q)=1.
   \end{equation}
\end{theorem}

Note that the number of distinct eigenvalues of $A$ is infinite and
dense in some interval!  However, the orders of the torsion subgroups
of $H$ is unbounded so that Corollary \ref{finiteev} is not
contradicted.  The eigenvalues of $A$ are algebraic numbers as
predicted by our Theorem 2.5 in \cite{MY}.  This can be seen as
follows: since $(\cos(p/q \pi)+i\sin(p/q\pi))^{q}=(-1)^p$, this
shows that  $\cos(p/q \pi)+i\sin(p/q\pi)$ is an algebraic number.
Therefore the real part, $\cos(p/q \pi)$, is also an algebraic number.

\section{The case of algebraic multipliers}
\label{sec:algmult}

The goal in this section is to extend the results that were obtained
in the previous sections, from rational multipliers to the more general
case of algebraic multipliers. Recall that a generic algebraic
multiplier is not
necessarily a rational multiplier. We start with an example of
algebraic numbers
on the unit circle that are not roots of unity. Consider the roots of
the polynomial,
\begin{equation}\label{eq:alg-rat}
   z^4 - z^3 + (2-k)z^2 - z +1
   =\left(z^2 - \frac{(1+\sqrt{4k+1})}{2}\cdot z + 1\right)\cdot
   \left(z^2 - \frac{(1-\sqrt{4k+1})}{2}\cdot z + 1\right)
\end{equation}
with $k$ a positive integer. This polynomial is irreducible over $\Z$
if $4k+1$ is not a square. We look for $k$ such that the first factor
has two distinct real roots while the second one has two complex
conjugate roots. Thus we seek $k$ so that $$
\frac{(1+\sqrt{4k+1})^2}{4} -4 > 0 \quad \mbox{and}
\quad \frac{(1-\sqrt{4k+1})^2}{4} -4 < 0.
$$
It is easily seen that the only values of $k$ satisfying these
conditions are $k=3,4,5$. For each of these choices,
two of the roots, denoted by $e^{i\theta}$ and $e^{-i\theta},$
lie on the unit circle and are roots of the second factor in
\eqref{eq:alg-rat}.
The two other roots are real,  denoted by $r$ and $r^{-1}$, where $r <
1$.
The numbers $e^{i\theta}$, $e^{-i\theta}$, $r$ and $r^{-1}$ are
algebraic integers, which are all
conjugate to each other.
Therefore $e^{i\theta}$ is not a root of unity since otherwise all its
conjugates
would also be roots of unity.
However, the numbers $e^{i\theta}$, $e^{-i\theta}$, $r$ and $r^{-1}$
are units in the corresponding ring of algebraic integers.
Since $e^{i\theta}$ is not a root of unity, its powers
are dense in the unit circle whereas the positive
powers of $r$ tend to 0.
For fixed $\alpha_1, \; \alpha_2 \in \mathbb R$ such that
$\theta=\alpha_2-\alpha_1$,
and for all $(m,n) \in \Z^2$, let
\begin{equation}
   \sigma((m',n'),(m,n))=\exp(-i(\alpha_1m'n+\alpha_2 n'm)).
\end{equation}
Then $\sigma$ is an algebraic multiplier on $\Z^2$
whose cohomology class $[\sigma] \in H^2(\Z^2, U(1)) \cong U(1)$
is equal to $e^{i\theta}$, so that $\sigma$ is {\em not} a rational
multiplier.
It is well known that $\sigma$ determines the noncommutative torus
$A_\theta$, see \cite{Boca}.

The trivial group has the $\sigma$-multiplier algebraic eigenvalue
property for any $\sigma$, since the eigenvalues are the zeros of the
characteristic polynomial. The same is true for every finite group. If
$G$ has the $\sigma$-multiplier algebraic eigenvalue property and $H$
is a subgroup of $G$, then $H$ also has the $\sigma$-multiplier
algebraic eigenvalue property.

\begin{theorem}
   \label{thm:free-maep}
   Every free group has the $\sigma$-multiplier algebraic eigenvalue
   property for every  $\sigma$.
\end{theorem}

\begin{proof}
   Let $G$ be a free group. Then $G$ has the algebraic eigenvalue
   property, corresponding to the identity multiplier, by
   Theorem 4.5 of \cite{dlmsy}. However for a free group, every
   multiplier is cohomologous to the identity, as free groups have
   no cohomology of degree two or higher. This can be seen by
   noting that the classifying space of a free group is a bouquet
   of circles, and so is one dimensional (see for example
   \cite[Chapter II, Section 4, Example 1]{Brown}.)
   By Corollary \ref{cor:cohom-maep} then, $G$ has the
   $\sigma$-multiplier algebraic eigenvalue for every algebraic
   multiplier $\sigma$.
\end{proof}

\begin{theorem}
   \label{thm:extension-maep}
   Suppose that we have a short exact sequence of groups
   \begin{equation}
     \label{eq:ses}
     1\to H\to G\stackrel{p}{\to} G/H\to 1
   \end{equation}
   where the quotient group $G/H$ is a finitely generated amenable   
   group.
   Let $\sigma'$ be an algebraic multiplier on $G/H$, and let
   $\sigma=p^\ast \sigma'$ be the pullback of $\sigma'$. Then
   if $H$ has the algebraic eigenvalue property,
   $G$ has the $\sigma$-multiplier algebraic eigenvalue property.
\end{theorem}

\begin{proof}
   We will show that an operator $A\in M_d(\Qbar(G,\sigma))$
   has only algebraic eigenvalues by first demonstrating that
   the point spectrum of $A$ is a subset of the union of the
   point spectra of a series $A^{(m)}$ of approximations to
   $A$, and then showing that each $A^{(m)}$ is equivalent
   to the untwisted action of a matrix over $\Qbar H$,
   and thus has only algebraic eigenvalues.
   
   In the following let $\calA$ be the von Neumann algebra
   $W_L^*(G,\sigma)\otimes M_d(\C)$, with trace $\tau=\tr_{G,\sigma}$
   as defined in Corollary~\ref{cor:trace}.  For finite $X\subset G/H$
   let $\calH_X$ be the subspace $l^2(p^{-1}(X))^d$ of $l^2(G)^d$
   and let $\calA_X$ be the commutant $B(\calH_X)^H$ of the right
   $H$-translations on $\calH_X$.
   Picking a right inverse $s$ of $p$, the isometry
   \begin{gather*}
      \iota_X: \calH_X=l^2(p^{-1}(X))^d \to l^2(H)^{d,\#X} \\
      (\iota_X f)(h)_{a,i} = f(s(x_i)h)_a\quad
      \text{for $h\in H$, $x_i\in X$, $a=1,\dots,d$}
   \end{gather*}
   induces an isomorphism $\psi_X$ from $\calA_X$ to
   $W_L^*(H)\otimes M_d(\C)\otimes M_{\#X}(\C)$.
   Define a trace $\tau_X$ on $A_X$ by
   \[
   \tau_X(B)=\frac{1}{\#X} (\tr_H\otimes\Tr) (\psi_X B)
   \]
   where $\tr_H$ is the usual trace on $W_L^*(H)$ and $\Tr$ is the
   canonical matrix trace on $M_d(\C)\otimes M_{\#X}(\C)$.
   In terms of the components $(B_{a,b})_{g,k}$ of an operator
   $B\in\calA_X$ (for $g,k\in G$ and $a,b=1,\dotsc,d$),
   the trace is given by
   \[
   \tau_X(B) =
   \frac{1}{\#X}\sum_{a=1}^d\sum_{x\in X} (B_{a,a})_{s(x),s(x)}.
   \]
   The von Neumann dimension associated with $\tau_X$ will be denoted
   by $\dim_X$.

   The multiplier $\sigma(g,h) = 1$ for all $h\in H$,
   so any operator $A\in\calA$ commutes
   with the right $H$-translations.
   For $A\in\calA$ and $X\subset G/H$ let $A^{(X)}=P_X A|_{\calH_X}$
   where $P_X$ is the orthogonal projection onto $\calH_X$.
   $A^{(X)}$ then belongs to $\calA_X$ and $\tau_X(A^{(X)})=\tau(A)$.

   The dimension functions on the algebras $\calA_X$ satisfy
   the following easily verifiable relations, for finite subsets
   $X\subseteq Y$:
   \begin{align}
      \label{eq:dimxincl}
      \dim_X L &\leq \dim_\tau N\quad
      \text{for all $L\eta\calA_X$, $N\eta\calA$ with $L\subset N$},
      \\
      \label{eq:dimxratio}
      \dim_X L &= \tfrac{\#Y}{\#X}\dim_Y L\quad
      \text{for all $L\eta \calA_X$},
      \\
      \label{eq:dimxproj}
      \dim_Y M&\geq \dim_Y P_X(M)\quad
      \text{for all $M\eta\calA_Y$}.
   \end{align}

   As $G/H$ is amenable and finitely generated, it admits
   a F\o{}lner exhaustion by finite subsets $\{X_m\}$
   such that
   \[
   \lim_{m\to\infty} \frac{\#\partial X_m}{\#X_m}=0,
   \]
   where $\partial X_m$ is the $\delta$-neighbourhood
   (with regard to the word metric on $G/H$) of $X_m$
   for any fixed $\delta$.

   In the following, let $A\in M_d(\Qbar(G,\sigma))$
   and let
   $A^{(m)}=A^{(X_m)}\in\calA_{X_m}$.
   Suppose $\lambda$ is not an eigenvalue of any of the
   $A^{(m)}$, and consider the space $E_\lambda$ of
   $\lambda$-eigenfunctions of $A$ with corresponding
   orthogonal projection $P_\lambda$. For any finite
   $X\subset G/H$,
   \begin{align*}
      \dim_\tau E_\lambda=\tau(P_\lambda)
      &=\tau_X(P_XP_\lambda|_{\calH_X})
      \\
      &\leq \dim_X \im P_X P_\lambda|_{\calH_X}\quad
      \text{( as $\norm{P_X P_\lambda}\leq 1$ ) }
      \\
      &\leq \dim_X P_X(E_\lambda).
   \end{align*}
   
   As $A$ is a matrix over the twisted group algebra, each
   component is a finite sum of twisted translations, and
   consequently $A$ has \emph{bounded propagation}.
   Explicitly, there are are only a finite number of the
   matrices of coefficients $A(g)\in M_d(\C)$ which are
   non-zero, and so we can choose a bound $\kappa$ by
   \[
   \kappa=\max \{d_{G/H}(1_{G/H},gH)\,|\ A(g)\neq 0\}
   \]
   (where $d_{G/H}$ is the word-metric on $G/H$) so that for
   $f$ with support in $p^{-1}(X)$, $Af$ will have support in
   $p^{-1}(X')$, where $X'=\{x\,|\ d_{G/H}(x,X)\leq \kappa\}$
   is the $\kappa$-neighbourhood of $X_m$.

   Let $X_m'$ be the $\kappa$-neighbourhood of $X_m$, and
   $\partial X_m=X_m' \setminus X_m$ be the outer
   $\kappa$-boundary of $X_m$. Then
   \[
   P_{X_m}A = P_{X_m}AP_{X_m'} =
   A^{(m)}P_{X_m} + P_{X_m}AP_{\partial X_m}.
   \]
   For $f\in E_\lambda$ with $P_{\partial X_m}f=0$ then,
   $P_{X_m} Af=\lambda P_{X_m} f = A^{(m)}P_{X_m}f$. By assumption
   though, $\lambda$ is not an eigenvalue of $A^{(m)}$, and so
   \begin{equation}
      \label{eq:zeroonbdy}
      f\in E_\lambda\,\text{ and } P_{\partial X_m}f=0
      \implies
      P_{X_m}f=0.
   \end{equation}
   Picking some superset $Y$ of $X_m'$, one has
   $P_{\partial X_m}P_Y=P_{\partial X_m}$ and $P_{X_m}P_Y=P_{X_m}$, and so
   (\ref{eq:zeroonbdy}) implies
   \begin{equation}
      \label{eq:kerinclusion}
      \ker P_{\partial X_m}|_{P_Y(E_\lambda)}
      \subseteq
      \ker P_{X_m}|_{P_Y(E_\lambda)}.
   \end{equation}
   Applying Lemma \ref{l:vndim2} gives
   \begin{align*}
      \dim_Y P_Y(E_\lambda)
      &=\dim_Y \ker P_{\partial X_m}|_{P_Y(E_\lambda)}
      +\dim_Y P_{\partial X_m}(E_\lambda)
      \\
      &=\dim_Y \ker P_{X_m}|_{P_Y(E_\lambda)}
      +\dim_Y P_{X_m}(E_\lambda),
   \end{align*}
   which with the inclusion (\ref{eq:kerinclusion}) in turn gives
   \[
   \dim_Y P_{X_m}(E_\lambda)\leq \dim_Y P_{\partial X_m}(E_\lambda)
   \leq \dim_Y \im P_{\partial X_m}=\frac{\#\partial X_m}{\#Y}.
   \]
   Then for any $m$ and $Y\supseteq X_m'$,
   \[
   \dim_\tau E_\lambda\leq \dim_{X_m}P_{X_m}(E_\lambda)
   =\frac{\#Y}{\#X_m}\dim_Y P_{X_m}(E_\lambda)
   \leq\frac{\#\partial X_m}{\#X_m},
   \]
   which goes to zero as $m$ goes to infinity, as the $X_m$ form
   a F\o{}lner exhaustion of $G/H$. Consequently
   $\dim_\tau E_\lambda=0$ and $\lambda$ is not an eigenvalue of
   $A$; that is, any eigenvalue of $A$ must be an eigenvalue of
   $A^{(m)}$ for some $m$.

   For $f\in \calH_X$, let $\hat f=\iota_X f$ be the
   corresponding element in $l^2(H)^{d,\#X}$
   defined by $\hat f(h)_{a,i}=f(s(x_i)h)_a$, as before.
   Then for every $x_i\in X$
   \begin{align*}
      (Af)(s(x_i)h)
      &=
      \sum_{j=1}^{\#X} \sum_{k\in H} A(s(x_i)hk^{-1}s(x_j)^{-1})
      f(s(x_i)k)\sigma(s(x_i)hk^{-1}s(x_j)^{-1},s(x_j)k)\\
      &=
      \sum_{j=1}^{\#X} \sum_{k\in H} A(s(x_i)hk^{-1}s(x_j)^{-1})
      f(s(x_i)k)\sigma'(x_i x_j^{-1},x_i)\\
      &=
      (B\hat f)(h)_i\in\C^d,
   \end{align*}
   where the matrix of coefficients of $h$ in $B\in M_{(d\#X)}(\Qbar H)$
   is given by
   \[
     B(h)_{(a,i),(b,j)}=
     A(s(x_i)hs(x_j)^{-1})_{a,b}\cdot\sigma'(x_i x_j^{-1},x_j)
   \]
   for $i,j=1,\dotsc,\#X$ and $a,b=1,\dotsc d$. As $H$ has the algebraic
   eigenvalue property, $B$ and hence $A^{X}$ have only algebraic  
   eigenvalues.
   Consequently, the operators $A^{(m)}$ and $A$ have
   only algebraic eigenvalues.
\end{proof}

One of the main theorems in this section is the following.

\begin{theorem}\label{fuchsian-maep}
   Let $\Gamma$ be the fundamental group of a closed Riemann surface of
   genus $g>1$.
   Then $\Gamma$ has the $\sigma$-multiplier algebraic eigenvalue
   property,
   where $\sigma$ is any algebraic multiplier on $\Gamma$.
\end{theorem}

We want to use Theorem \ref{thm:extension-maep} using the exact
sequence of Remark (\ref{surface}). To do this, it is necessary to
prove that every algebraic multiplier $\sigma$ on $\Gamma$ is
cohomologous to the pull-back of an algebraic multiplier $\sigma'$ on
$\mathbb Z^{2g}$.  The construction of $\sigma'$ was used in
\cite[Section 7.2]{CHMM98}.  We follow it closely paying particular
attention to algebraicity.

Recall that the {\em area cocycle} $c$ of the fundamental group
of a compact Riemann surface,
$\Gamma=\Gamma_g$ is a canonically defined 2-cocycle
on $\Gamma$ that is defined as follows.
Firstly, recall the definition of
a well known area 2-cocycle on ${\bf PSL}(2, \R)$.
${\bf PSL}(2, \R)$ acts on $\mathbb H$ so that
$\mathbb H \cong {\bf PSL}(2, \R)/{\bf SO}(2)$.
Then
$$c(\gamma_1, \gamma_2) = \text{Area}_{\mathbb H}(\Delta(o,
\gamma_1\cdot o,
{\gamma_1\gamma_2}\cdot o)),$$
where $o$
denotes an origin in $\mathbb H$ and
$\text{Area}_{\mathbb H}(\Delta(a,b,c))$ is the oriented hyperbolic
area of the
geodesic triangle in $\mathbb H$
with vertices at $a, b, c \in \mathbb H$. The restriction of
$c$ to the subgroup $\Gamma$ is the  area cocycle $c$ of $\Gamma$.
We use the additive notation when discussing area cocycles and remark
that $(2\pi)^{-1}c$ represents an integral class in $H^2(\Gamma,
{\mathbb
   R}) \cong {\mathbb R}$ as follows from Gauss-Bonnet theorem.

Let $\Omega_j$ denote the (diagonal)
operator on $l^2(\Gamma)$ defined by
$$
\Omega_j f(\gamma) = \Omega_j(\gamma)f(\gamma) \quad \forall f\in
l^2(\Gamma)
\quad \forall \gamma\in \Gamma
$$
where
$$
\Omega_j(\gamma) = \int_o^{\gamma \cdot o} \alpha_j \quad j=1,\ldots ,2g
$$
and where
\begin{equation} \label{sympl:basis}
   \{\alpha_j \}_{j=1}^{2g}=\{ a_j \}_{j=1} ^g \cup \{ b_{j}
   \}_{j=1}^g \end{equation}
is a collection of harmonic $1$-forms on the compact Riemann surface
$\Sigma_g = \hyp/\Gamma$, generating $H^1(\Sigma_g,\R)=\R^{2g}$.
We abuse the notation slightly and do not distinguish between a form on
$\Sigma_g$ and its pullback to the hyperbolic plane as well as between
an element of $\Gamma$ and a loop in $\Sigma_g$ representing it.

Notice that we can write equivalently
$$
\Omega_j(\gamma) = c_j(\gamma),
$$
where the group cocycles $c_j$ form a symplectic basis for
$H^1(\Gamma,\Z)=\Z^{2g}$, with generators $\{ \alpha_j
\}_{j=1,\ldots ,2g}$, as in (\ref{sympl:basis}) and can be defined as
the integration on loops on $\Sigma_g$,
$$
c_j(\gamma)=\int_\gamma \alpha_j.
$$

Define
$$
\Psi_j(\gamma_1, \gamma_2) = \Omega_j(\gamma_1) \Omega_{j+g}(\gamma_2)
- \Omega_{j+g}(\gamma_1) \Omega_j(\gamma_2).
$$

Let $\Xi: \mathbb H \to {\mathbb R}^{2g}$ denote the Abel-Jacobi
map
$$ \Xi: x \mapsto \left( \int_o^x a_1, \int_o^x
   b_{1}, \ldots, \int_o^x a_g, \int_o^x b_{g} \right), $$
where $\displaystyle\int_o^x$ means integration along the unique
geodesic in
${\mathbb H}$ connecting $o$ to $x$. Having chosen an origin $o$ once
and for all we make an identification $\Gamma \cdot o \cong\Gamma$.
Note that $\Gamma$ acts on ${\mathbb R}^{2g}$ in a natural way and 
the map $\Xi$ is $\Gamma$-equivariant. In addition,
the map $\Xi$ is a symplectic map, that is, if $\omega$ and $\omega_J$
are the respective symplectic
2-forms, then one has $\Xi^*(\omega_J) = k\omega$ for a suitable
constant $k$. Henceforth, we renormalize $\omega$ (and consequently the
area cocycle $c$) so that $\Xi^*(\omega_J) = \omega$
One then has the following
geometric lemma \cite{CHMM98}, \cite{MM01}.

\begin{lemma}
   \label{geom}
   $$
   \Psi (\gamma_1, \gamma_2) = \sum_{j=1}^g  \Psi_j(\gamma_1, \gamma_2)
   = \int_{\Delta_E(\gamma_1, \gamma_2)} \omega_J
   $$ where $\Delta_E(\gamma_1, \gamma_2)$ denotes the Euclidean
   triangle with vertices at $\Xi(o), \Xi(\gamma_1\cdot o)$ and $
   \Xi(\gamma_1\gamma_2\cdot o)$, and $\omega_J$ denotes the flat
   K\"ahler 2-form on the universal cover of the Jacobi variety. That
   is, $\sum_{j=1}^g \Psi_j(\gamma_1, \gamma_2)$ is equal to the
   Euclidean area of the Euclidean triangle
   $\Delta_E(\gamma_1,\gamma_2)$.
\end{lemma}

That is, the cocycle $\Psi = p^*(\Psi')$, where $\Psi'$ is a 2-cocycle
on $ \mathbb Z^{2g}$ and $p$ is defined as the projection.
$$
1\to F\to \Gamma \stackrel{p}{\to} \mathbb Z^{2g}\to 1
$$

The following lemma is also implicit in \cite{CHMM98}, \cite{MM01}.

\begin{lemma}
   \label{comparison}
   The hyperbolic area group $2$-cocycle $c$ and
   the Euclidean area group $2$-cocycle $\Psi$
   on $\Gamma$, are cohomologous.
\end{lemma}

\begin{proof}

   Observe that since $\omega = \Xi^* \omega_J$, one has
   $$
   c(\gamma_1, \gamma_2) = \int_{\Delta(\gamma_1, \gamma_2)}
   \omega = \int_{\Xi(\Delta(\gamma_1, \gamma_2))}
   \omega_J.
   $$
   Therefore the difference
   \begin{align*}
     \Psi (\gamma_1, \gamma_2)   - c(\gamma_1, \gamma_2)
     & = \int_{\Delta_E(\gamma_1, \gamma_2)}
     \omega_J  - \int_{\Xi(\Delta(\gamma_1, \gamma_2))}
     \omega_J\\
     & = \int_{\partial \Delta_E(\gamma_1, \gamma_2)}
     \Theta_J  - \int_{\partial\Xi(\Delta(\gamma_1, \gamma_2))}
     \Theta_J,
   \end{align*}
   where $\Theta_J$ is a 1-form on the universal cover ${\mathbb R}^{2g}$
   of the
   Jacobi variety such that $d\Theta_J = \omega_J$.

   Let
   $h(\gamma) = \int_{\Xi(\ell (\gamma))} \Theta_J -
   \int_{m(\gamma)}\Theta_J $,
   where $\ell(\gamma)$ denotes the unique geodesic in
   $\mathbb H$ joining $ o$ and $\gamma\cdot o$ and $m(\gamma)$ is the
   straight line in the Jacobi variety joining the points $\Xi(o)$ and
   $\Xi(\gamma\cdot o)$.  We can also write
   $h(\gamma)=\int_{D(\gamma)} \omega_J$, where $D(\gamma)$ is an
   arbitrary
   topological disk in
   ${\mathbb R}^{2g}$ with boundary $\Xi(\ell(\gamma))\cup m(\gamma)$.
   Thus the equality above can be rewritten as

\begin{align*}
     \Psi (\gamma_1, \gamma_2)  - c(\gamma_1, \gamma_2)
     & = h(\gamma_1)  - h(\gamma_1\gamma_2)  +
     \int_{\gamma_1\cdot\Xi(l(\gamma_2))}\Theta_J - 
     \int_{\gamma_1\cdot m(\gamma_2)}\Theta_J\\
     & = h(\gamma_1)  - h(\gamma_1\gamma_2)  +
     \int_{D(\gamma_2)}(\gamma_1)^*d\Theta_J \\
     & = h(\gamma_1)  - h(\gamma_1\gamma_2)  +
     h(\gamma_2) \\
     & = \delta h (\gamma_1, \gamma_2)
   \end{align*}
since $\omega_J = d \Theta_J$ is invariant under the action of $\Gamma$.

\end{proof}

\begin{lemma}
   \label{comparison2}
   Let $\Gamma$ be the fundamental group of a closed, genus $g$ Riemann
   surface
   and
   $$
   1\to F\to \Gamma \stackrel{p}{\to} \mathbb Z^{2g}\to 1
   $$
   where $F$ is a free group  and $ \mathbb Z^{2g}$ the free
   abelian group  as in Remark \ref{surface}. Then
   every multiplier $\sigma$ on $\Gamma$ is cohomologous
   to a multiplier $\sigma' = p^* (\sigma'')$ on $\Gamma$
   where $\sigma''$ is a multiplier on $ \mathbb Z^{2g}$.

   In addition,
   every algebraic multiplier $\sigma$ on $\Gamma$ is cohomologous
   to an algebraic multiplier $\sigma' = p^* (\sigma'')$ on $\Gamma$
   where $\sigma''$ is an algebraic multiplier on $ \mathbb Z^{2g}$.
   More precisely, if  $\sigma \in
   Z^2(\Gamma, \mathcal U(\Qbar))$ then $\sigma''$ can be chosen from
   $Z^2(\mathbb Z^{2g}, \mathcal U(\Qbar))$ so that $\sigma$ and $\sigma'$ are
   cohomologous in $Z^2(\Gamma,\mathcal U(\Qbar) )$.
\end{lemma}

\begin{proof}
   Observe that $\Xi(\Gamma\cdot o )\subset {\mathbb Z}^{2g}\subset
   {\mathbb R}^{2g}$. It follows that the Euclidean area cocycle and
   its pullback represent integral cohomology classes. By the lemma
   above, the cohomology class of $c$ is integral. Now let $\sigma$ be
   an arbitrary multiplier on $\Gamma$. Since $H^2(\Gamma,
   A)=A$ for every abelian group $A$ and $H^3(\Gamma, {\mathbb Z}) =
   0$, we see that $\sigma$ is cohomologous to a multiplier
   $\sigma_1 = \exp(2\pi i \theta c)$, where $\theta$ is a real number.
   By Lemma \ref{comparison} we see that $\sigma_1$ is
   cohomologous to $ p^* (\sigma'')$, where
   $\sigma'' =  \exp(2\pi i \theta\Psi')$
   is a multiplier on $ \mathbb Z^{2g}$. 
   
   To prove the last claim, we
   identify the the group cohomology with the cohomology of the surface
   $\Sigma_g={\mathbb H}/\Gamma$. Since the value of the cocycle on the
   fundamental class depends only on the cohomology class and
   $c(\Sigma_g)=2g-2$, we see that $\sigma'(\Sigma_g) = \exp(2\pi
   i\theta)^{2g-2} =\sigma(\Sigma_g)$ is algebraic. It follows that
   $\exp(2\pi i \theta)$ is an algebraic number so that $\sigma''$ is an
   algebraic cocycle. Now both $\sigma$ and $\sigma'$ are algebraic
   cocycles. They are cohomologous in $Z^2(\Gamma, U(1))$. For
   any coefficients the cohomology class of the cocycle is
   determined by the value of the cocycle on the fundamental class.
   Therefore 
   $\sigma$ and $\sigma'$ represent the same cohomology class in
   $H^2(\Gamma,\mathcal U(\Qbar))$ i.e.\ are cohomologous in
   $Z^2(\Gamma,\mathcal U(\Qbar))$.
\end{proof}

\begin{proof}[{\bf Proof of Theorem \ref{fuchsian-maep}}]
   Recall that if two multipliers $\sigma$, $\sigma'$ on $\Gamma$ are
   cohomologous, then $\Gamma$ has the $\sigma$-multiplier algebraic
   eigenvalue property if and only if $\Gamma$ has the
   $\sigma'$-multiplier algebraic eigenvalue property (Corollary
   \ref{cor:cohom-maep}.)  Since the free group $F$ has the algebraic
   eigenvalue property, and since $ \mathbb Z^{2g}$ is a finitely
   generated amenable group, by applying Theorem
   \ref{thm:extension-maep} and Lemma \ref{comparison2}, we deduce
   Theorem \ref{fuchsian-maep}.
\end{proof}

\section{Generalized integrated density of states and spectral gaps}

In this section, we will realize the von Neumann trace on the group von
Neumann algebra
of a surface group, as a generalized  integrated density of states,
which is an important step to relating it
directly to the physics of the quantum Hall effect.

Our first main theorem is the following.

\begin{theorem}
   \label{thm:ids-extension}
   Consider the situation of Theorem \ref{thm:extension-maep},
   where we have a short exact sequence of groups
   \begin{equation}
     \label{eq:ids-ses}
     1\to H\to G\stackrel{p}{\to} G/H\to 1
   \end{equation}
   where the quotient group $G/H$ is finitely generated and amenable.
   Let $\sigma'$ be a  multiplier on $G/H$, and let
   $\sigma=p^\ast \sigma'$ be the pullback of $\sigma'$.
   Let $A\in M_d(\C(G,\sigma))$ be a self-adjoint operator acting
   on $l^2(G)^d$, being a member of the von Neumann algebra
   $\calA=W_L^*(G,\sigma)\otimes M_d(\C)$ with trace
   $\tau=\tr_{G,\sigma}$.
   
   For finite subsets $X$ of $G/H$, let $\calH_X=l^2(p^{-1}(X))^d$
   be the space of functions with support on $p^{-1}(X)$,
   and let $\calA_X=B(\calH_X)^H$ be the commutant of the right
   $H$-translations on $\calH_X$. Pick a right inverse $s$ of the
   projection $p$ and give $\calA_X$ the trace $\tau_X$
   as in the proof of Theorem \ref{thm:extension-maep}, which
   in terms of the components $(B_{a,b})_{g,k}$ of an operator
   $B\in\calA_X$ (for $g,k\in G$ and $a,b=1,\dotsc,d$)
   is given by
   \[
   \tau_X(B) =
   \frac{1}{\#X}\sum_{a=1}^d\sum_{x\in X} (B_{a,a})_{s(x),s(x)}.
   \]
   Let $A^{(X)}=P_X A|_{\calH_X}\in\calA_X$ where $P_X$ is the
   orthogonal projection onto $\calH_X$.   
   
   Choose a F\o{}lner exhaustion ${X_m}$ of $G/H$. Then the spectral
   density function of $A$ equals the generalised integrated density
   of states as given by the (normalised) spectral density functions of
   the operators $A^{(m)}=A^{(X_m)}$. That is, with spectral density
   functions $F$ of $A$ and $F_m$ of the $A_m$,
   \begin{align*}
   F_m(\lambda) &= \tau_{X_m} (\chi_{(-\infty,\lambda]}(A^{(m)})),
   &
   F(\lambda) &= \tau (\chi_{(-\infty,\lambda]}(A)),
   \end{align*}
   the $F_m$ converge point-wise to $F$ at every $\lambda$,
   \begin{equation}
      \label{eq:gids}
      \lim_{m\to\infty} F_m(\lambda)=
      F(\lambda)\qquad\forall\lambda\in\R.
   \end{equation}
\end{theorem}

The proof of this theorem in the case of $H=1$ was given in \cite{MY}
and \cite{MSY} for the discrete magnetic Laplacian. To establish
this theorem in our more general situation, we apply the same arguments,
slightly generalized as follows, relying upon the notation established
in the proof of Theorem \ref{thm:extension-maep}.

\begin{lemma}
   For any polynomial $p$
   \[
   \lim_{m\to\infty} \tau_m(p(A^{(m)}))=\tau(p(A)).
   \]
\end{lemma}
\begin{proof}
   The argument is exactly that of \cite{MY}, Lemma~2.1, and
   relies upon the amenability of $G/H$.
\end{proof}

\begin{lemma}
   \label{lem:wkspecapprox}
   Suppose $f(\lambda)$ and $f_m(\lambda)$ ($m=1,2,\dotsc$) are
   monotonically increasing right continuous functions on $\R$ that
   are zero for $\lambda<a$ and constant for $\lambda\geq b$, for fixed
   $a$ and $b$.
   Further suppose that
   \begin{equation}
      \label{eq:polylimeq}
      \lim_{m\to\infty} \int p\ df_m = \int p\ df
   \end{equation}
   for all polynomials $p$, where the integrals are Lebesgue-Stieltjes
   integrals.  Then
   \begin{equation*}
      f(\lambda)
      =\finfplus(\lambda)=\fsupplus(\lambda)
   \end{equation*}
   for all $\lambda$, where $\fsup^+$ and $\finf^+$ are defined in terms
   of the $f_m$ by
   \begin{equation}
      \label{eq:finfdef}
      \begin{aligned}
         \finf(\lambda)
         &=\liminf_m f_m(\lambda),
         &
         \finfplus(\lambda)
         &=\lim_{\epsilon\to 0^+} \finf(\lambda+\epsilon),
         \\
         \fsup(\lambda)
         &=\limsup_m f_m(\lambda),
         &
         \fsupplus(\lambda)
         &=\lim_{\epsilon\to 0^+} \fsup(\lambda+\epsilon).
      \end{aligned}
   \end{equation}
   In particular $f(\lambda)=\lim_{m\to\infty}
   f_m(\lambda)$ at all points of continuity of $f$, which is at all
   but a countable number of points.
\end{lemma}

\begin{proof}
   The proof follows that of part (i) of Theorem 2.6 of
   \cite{MY}.
   
   Take a sequence of successively closer polynomial approximations
   $p_j$ to the characteristic function $\chi_{(-\infty,x]}$ over the
   interval $[a,b)$ such that
   \[
   \chi_{(-\infty,x]}(\lambda)\leq p_j(\lambda)
   \leq
   \chi_{(-\infty,x+\tfrac{1}{j}]}(\lambda)+\tfrac{1}{j}
   \quad\forall \lambda\in {[a,b)},\ j\geq 1.
   \]
   Then for all $j$,
   \begin{gather}
      \label{eq:pjdfm}
      f_m(x)\leq
      \int_a^b p_j(\lambda)df_m(\lambda)
      \leq f_m(x+\tfrac{1}{j})+\tfrac{1}{j}(b-a),
      \\
      \label{eq:pjdf}
      f(x)\leq
      \int_a^b p_j(\lambda)df(\lambda)
      \leq f(x+\tfrac{1}{j})+\tfrac{1}{j}(b-a).
   \end{gather}
   Taking the limit as $m$ goes to infinity, equations (\ref{eq:pjdfm})
   and (\ref{eq:polylimeq}),
   \begin{equation}
      \label{eq:flim1}
      \fsup(x)\leq \int_a^b p_j(\lambda)df(\lambda)
      \leq \finf(x+\tfrac{1}{j})+\tfrac{1}{j}(b-a)
      \quad\forall j\geq 1.
   \end{equation}
   The right continuity of $f$ with equation (\ref{eq:pjdf}) gives
   \[
   \lim_{j\to\infty} \int_a^b p_j(\lambda)df(\lambda)
   = f(x),
   \]
   and so taking the limit of (\ref{eq:flim1}) as $j$ goes to infinity
   gives
   \begin{equation}
      \label{eq:flim2}
      \fsup(x)\leq f(x)\leq \finfplus(x)\quad\forall x.
   \end{equation}
   Again using the right continuity of $f$,
   \[
   f(x)\leq\finfplus(x)\leq\fsupplus(x)\leq\fplus(x)=f(x).
   \]
   $f(x)$ is monotonically increasing in $x$ and bounded, so can
   have at most a countable number of discontinuities. If $f$ is
   continuous at $x$ then equation (\ref{eq:flim2}) implies that
   $\finf(x)=f(x)=\fsup(x)$.
\end{proof}

\OMIT{
   When $f$ and $f_m$ are spectral density functions, the situation
   of the preceding lemma is called \emph{weak spectral convergence}.
}

\begin{lemma}
   \label{lem:Fwkspecapprox}
   Let $F$ and $F_m$ be as in the statement of
   Theorem~\ref{thm:ids-extension}. Then using the notation
   (\ref{eq:finfdef}) of Lemma~\ref{lem:wkspecapprox},
   \[
   F(\lambda)=\Fsupplus(\lambda)=\Finfplus(\lambda)
   \quad\forall\lambda\in\R,
   \]
   with
   \begin{equation}
      \label{eq:wklim}
      \lim_{m\to\infty} F_m(\lambda)=F(\lambda)
      \qquad\text{$\forall \lambda\in\R$
         such that F is continuous at $\lambda$.}
   \end{equation}
\end{lemma}
\begin{proof}
   This is an immediate consequence of the two preceding lemmas.
\end{proof}

The convergence (\ref{eq:wklim}) can be extended to all $\lambda$
by showing that the jumps of the spectral density functions 
at points of discontinuity also converge.

\begin{lemma}[Corollary 3.2 of \cite{MSY}]
   \label{lem:onjumpconv}
   Let $f$ and $f_m$ (for $m=1,2,\dotsc$) be monotonically increasing
   right continuous functions on $\R$ satisfying
   $f(\lambda)=\finfplus(\lambda)=\fsupplus(\lambda)$ at all $\lambda$,
   as in Lemma~\ref{lem:wkspecapprox}. Denote the jumps at $\lambda$
   of $f$ and the $f_m$ by $j$ and $j_m$ respectively,
   \begin{align*}
      j_m(\lambda)
      &=\lim_{\epsilon\to 0^+} f_m(\lambda)-f_m(\lambda-\epsilon),
      \\
      j(\lambda)
      &=\lim_{\epsilon\to 0^+} f(\lambda)-f(\lambda-\epsilon).
   \end{align*}
   Suppose the $j_m$ converge to $j$ point-wise at all $\lambda$.
   Then the $f_m$ converge to $f$ point-wise at all $\lambda$.
\end{lemma}

To obtain point-wise convergence of $f_m$ to $f$ at every point,
it is in fact sufficient to show that
$\liminf_m j_m(\lambda)\geq j(\lambda)$ at all $\lambda$, due to
the following lemma.

\begin{lemma}
   \label{lem:jmupperbound}
   Let $f$ and $f_m$ (for $m=1,2,\dotsc$) be monotonically increasing
   right continuous functions on $\R$ satisfying
   $f(\lambda)=\finfplus(\lambda)=\fsupplus(\lambda)$ at all $\lambda$,
   as in Lemma~\ref{lem:wkspecapprox}. Denote the jumps at $\lambda$
   of $f$ and $f_m$ by $j$ and $j_m$ respectively, as in
   Lemma~\ref{lem:onjumpconv}. Then
   \[
   \limsup_m j_m(\lambda)\leq j(\lambda)\quad\forall\lambda\in\R.
   \]
\end{lemma}
\begin{proof}
   Fix $\lambda$. By monotonicity,
   \begin{equation}
      \label{eq:jm1}
      j_m(\lambda)\leq f_m(\lambda+\epsilon)-f_m(\lambda-\epsilon)
      \quad\forall\epsilon>0.
   \end{equation}
   $f$ is continuous at all but a countable number of points,
   and at points $x$ of continuity, $f_m(x)\to f(x)$ as
   $m\to\infty$. Pick a decreasing sequence $\epsilon_k\to 0$ such
   that $f$ is continuous at $\lambda+\epsilon_k$ and
   $\lambda-\epsilon_k$ for all $k$. Then taking the limit in $m$
   of (\ref{eq:jm1}) gives
   \[
   \limsup_m j_m(\lambda)\leq
   f(\lambda+\epsilon_k)-f(\lambda-\epsilon_k)\quad\forall k.
   \]
   By right continuity of $f$,
   $f(\lambda+\epsilon_k)-f(\lambda-\epsilon_k)$ converges to
   $j(\lambda)$ from above as $k$ goes to infinity. Thence
   on taking the limit in $k$,
   $\limsup_m j_m(\lambda)\leq j(\lambda)$.
\end{proof}

Now consider the situation of Theorem~\ref{thm:ids-extension}.
We already have a weak spectral approximation by virtue of
Lemma~\ref{lem:Fwkspecapprox}, so all we require now is to show
convergence of the jumps in $F_m$ to those of $F$.

\begin{theorem}
   \label{thm:Djumpconv}
   Let $D(\lambda)$ and $D_m(\lambda)$ denote the jumps at $\lambda$
   of the spectral density functions $F$ and $F_m$ respectively.
   Then
   \[
   \lim_{m\to\infty} D_m(\lambda) = D(\lambda)\quad\forall\lambda\in\R.
   \]
\end{theorem}
\begin{proof}
   Let $\dim_X$ be the von Neumann dimension associated with
   the trace $\tau_X$ on $\calA_X$.
   Note that $\dim \calH = \dim \ker B + \dim \im B$ for an operator
   $B$ in a von Neumann algebra of operators acting on a Hilbert space
   $\calH$, with finite von Neumann dimension $\dim$. So
   \begin{gather*}
      D(\lambda)=
      \dim_\tau \ker (A-\lambda)=
      d- \dim_\tau \im (A-\lambda),
      \\
      D_m(\lambda)=
      \dim_{X_m} \ker (A^{(m)}-\lambda) =
      d-\dim_{X_m} \im (A^{(m)}-\lambda).
   \end{gather*}
   As in the proof of Theorem~\ref{thm:extension-maep}, let $\kappa$
   be the propagation of the operator $A$ with respect to the word
   metric $d_{G/H}$ on $G/H$ and let $X_m'$ be the $\kappa$-neighbourhood
   of $X_m$ so that $f\in\calH_{X_m}$ implies $Af\in\calH_{X_m'}$;
   equivalently, $P_{X_m'}A|_{\calH_{X_m}}=A|_{\calH_{X_m}}$ for all $m$.

   The space $\im (A-\lambda)|_{\calH_{X_m}}$ is
   affiliated with $\calA_{X_m'}$. Recall the properties
   (\ref{eq:dimxincl}), (\ref{eq:dimxratio}) and (\ref{eq:dimxproj})
   of $\dim_X$ as listed in the proof of Theorem~\ref{thm:extension-maep}.
   Then
   \begin{align*}
      \dim_{X_m'} \im (A^{(m)}-\lambda)
      &=
      \dim_{X_m'} P_{X_m}(\im (A-\lambda)|_{\calH_{X_m}})\\
      &\leq
      \dim_{X_m'} \im (A-\lambda)|_{\calH_{X_m}}\\
      &\leq
      \dim_\tau \im (A-\lambda),
   \end{align*}
   and
   \[
   \dim_{X_m}\im (A^{(m)}-\lambda)=
   \frac{\#X_m'}{\#X_m}\dim_{X_m'}\im (A^{(m)}-\lambda).
   \]
   The $X_m$ constitute a F\o{}lner exhaustion of $G/H$ and so
   $\frac{\#X_m'}{\#X_m}$ tends to $1$ as $m$ goes to infinity.
   Taking limits gives
   \begin{align*}
      \liminf_m D_m(\lambda)
      &=d-\limsup_m\ \dim_{X_m}\im(A^{(m)}-\lambda)
      \\
      &\geq d-\dim_\tau\im(A-\lambda)
      =D(\lambda).
   \end{align*}
   Finally, applying Lemma~\ref{lem:jmupperbound} gives
   \[
   D(\lambda)
   \leq \liminf_m D_m(\lambda)
   \leq \limsup_m D_m(\lambda)
   \leq D(\lambda).
   \]
\end{proof}

The proof of Theorem~\ref{thm:ids-extension} now follows from
Lemmas~\ref{lem:Fwkspecapprox} and \ref{lem:onjumpconv},
and Theorem~\ref{thm:Djumpconv}.

The following corollary is an immediate consequence of
Theorem~\ref{comparison2} and Theorem \ref{thm:ids-extension}.

\begin{corollary}[Generalized IDS]
   \label{fuchsian approx}
   Let $G= \Gamma$ be the fundamental group of a closed Riemann
   surface of genus $g>1$, $G/H = \Z^{2g}$ be the abelianisation of
   $G$, in which case the commutator subgroup $H=F$ is a free group.
   Then the equality between the generalized integrated density of
   states and the von Neumann spectral density function given in
   equation \eqref{eq:gids} holds for every multiplier $\sigma$
   on $\Gamma$.
\end{corollary}

\begin{corollary}[Criterion for spectral gaps]
   \label{spectral gap}
   Consider the situation in Theorem \ref{thm:ids-extension}. The
   interval $(\lambda_1, \lambda_2)$ is in a gap in the spectrum of
   $A$ if and only if
   \begin{equation}\label{thm: spectral gaps}
     \lim_{m\to\infty} \left(F_m(\lambda_2) - F_m(\lambda_1)\right) = 0.
   \end{equation}
\end{corollary}

\begin{proof}
   The interval $(\lambda_1, \lambda_2)$ is in a gap in the spectrum
   of $A$ if and only if $F(\lambda_2) = F(\lambda_1)$. By Theorem
   \ref{thm:ids-extension}, this is true if and only if
   $$
   \lim_{m\to\infty} \left(F_m(\lambda_2) - F_m(\lambda_1)\right)  =
   F(\lambda_2) - F(\lambda_1) = 0.
   $$
\end{proof}

\section{The class $\calK$ and extensions with cyclic kernel}
\label{sec:cyclic}
\label{sec:classK}

In this section we prove the results cited in the earlier sections
concerning the class of groups
$\calK$. Namely we show that the class $\calK$ is closed under taking
extensions with cyclic kernel, and that every group in
$\calK$ has the algebraic eigenvalue property.

Recall that the class $\calK$ is the smallest class of groups
which contains the free groups and the amenable groups, and is
closed under directed unions and under
taking extensions with amenable quotients. It can be seen that
every group in $\calK$ must belong to some $\calK_\alpha$
defined inductively as follows.

\begin{definition}
   Define the nested classes $\calK_\alpha$, $\alpha$ an ordinal, by
   \begin{itemize}
   \item
     $\calK_0$ consists of all free groups and all
     discrete amenable groups,
   \item
     $\calK_{\alpha+1}$ consists of all extensions of groups
     in $\calK_\alpha$ with amenable quotient, and all directed unions
     of groups in $\calK_\alpha$,
   \item
     $\calK_\beta=\bigcup_{\alpha<\beta} \calK_\alpha$ when
     $\beta$ is a limit ordinal.
   \end{itemize}
   A group is in $\calK$ if and only if it is in a class $\calK_\alpha$
   for some ordinal $\alpha$.
\end{definition}

Recall Theorem \ref{thm:Kalg}:

\begin{theorem*}
   Every group in $\calK$ has the algebraic eigenvalue property.
\end{theorem*}

The proof follows closely that of Corollary 4.8 of \cite{dlmsy},
and relies upon the same key theorem:

\begin{theorem*}[4.7 of \cite{dlmsy}]
   Let $H$ have the algebraic eigenvalue property. Let $G$
   be a generalized amenable extension of $H$. Then $G$ has
   the algebraic eigenvalue property.
\end{theorem*}

The proof then proceeds by transfinite induction.

\begin{proof}[Proof of Theorem \ref{thm:Kalg}]
   The algebraic eigenvalue property holds for groups in $\calK_0$ by
   Corollary 4.8 in \cite{dlmsy}.

   Proceeding by transfinite induction, suppose $\calK_\alpha$ has
   the algebraic eigenvalue property for all $\alpha$ less than some
   given ordinal $\beta$.

   When $\beta=\alpha+1$ for some ordinal $\alpha$, a group $G$ is in
   $\calK_\beta$ if it is an extension of a group in $\calK_\alpha$ with
   amenable quotient, or is the directed union of groups $G_i \in
   \calK_\alpha$.

   Note that $A\in M_d(\Qbar G)$ can be regarded
   as a matrix $A'$ in $M_d(\Qbar H)$ where $H$
   is a finitely generated subgroup of $G$, generated by the finite
   support of the $A_{i,j}$ in $G$.  By Proposition 3.1 of
   \cite{Schick(1998a)}, the spectral density functions of $A$ and $A'$
   coincide. As subgroups of a group with the algebraic eigenvalue
   property also have the property, it follows that a group has the
   algebraic eigenvalue property if and only if it holds for all of its
   finitely generated subgroups.  If $G\in\calK_{\alpha+1}$ is the
   directed union of groups $G_i\in \calK_\alpha$, it follows that
   every finitely generated subgroup of $G$ is in some $G_i$ and so has
   the algebraic eigenvalue property. $G$ therefore has the algebraic
   eigenvalue property.

   Suppose that $G$ is the extension of a group $H$ in $\calK_\alpha$
   with amenable quotient $R=G/H$, i.e. $H\to G\stackrel{p}{\to} R$,
   then $R = \bigcup_{j\in J} R_j$  is a directed union,
   where $R_j, \; j\in J$ are finitely
   generated and amenable (since amenability is subgroup closed).
   Consider the extensions $G_j = p^{-1}(R_j)$ of $H$, with finitely
   generated amenable quotient: then Theorem 4.12 of \cite{dlmsy}
   applies to show that $G_j$ has the algebraic eigenvalue property.

   But $G= \bigcup_{j\in J} G_j$ is the directed union
   of groups $G_j$ having the algebraic eigenvalue property,
   so $G$ also has the algebraic
   eigenvalue property by the argument given in the previous paragraph.

   Finally, let $\beta$ be a limit ordinal, so that
   $\calK_\beta=\bigcup_{\alpha<\beta} \calK_\alpha$.
   If $G$ is in $\calK_\beta$, then it is in $\calK_\alpha$ for some
   $\alpha<\beta$, and so has the algebraic eigenvalue property by
   the induction hypothesis.

   Therefore groups in $\calK_\beta$ have the algebraic eigenvalue
   property,
   and the result follows by induction.
\end{proof}


We next address Proposition \ref{prop:cyclic}.
Its proof follows from the following lemmas.

\begin{lemma}\label{subgroup}
   The class $\calK_\alpha$ is subgroup-closed for all ordinals
   $\alpha$.
\end{lemma}

\begin{proof}
   The classes of free groups and amenable groups are both closed
   under taking subgroups, so $\calK_0$ is subgroup-closed.

   Suppose $\calK_\alpha$ is subgroup closed for all $\alpha<\beta$,
   let $G\in\calK_\beta$ and let $H$ be a subgroup of $G$.
   If $\beta$ is a limit ordinal then $G\in\alpha$ for some
   $\alpha<\beta$, and so $H\in\calK_\alpha$, and so in $\calK_\beta$.

   Suppose then that $\beta$ is not a limit ordinal, with
   $\beta=\alpha+1$ for some ordinal $\alpha$.
   If $G$ is a directed union of groups $G_j$ in $\calK_\alpha$,
   then $H$ is a directed union of groups $H\cap G_j$ which are also in
   $\calK_\alpha$ by the induction hypothesis.

   $G$ must otherwise be an extension of a group $N$ in $\calK_\alpha$
   with amenable quotient $B$. In this case, $H$ is an extension of
   $N\cap H$ with quotient $H/(N\cap H)\cong NH/N \subset G/N \cong B$,
   and hence is an extension with amenable quotient of $N\cap H$ which
   is in $\calK_\alpha$ by the hypothesis. $H$ is therefore in
   $\calK_\beta$.
\end{proof}

\begin{lemma}
   Let $G$ be a group whose $j$-th derived group $\derived[j]{G}$
   is in $\calK_\alpha$. Then $G\in\calK_{\alpha+j}$.
\end{lemma}

\begin{proof}
   $G$ is an extension with kernel its derived group $\derived{G}$
   by the short exact sequence
   $\derived{G}\longrightarrow G\longrightarrow G/\derived{G}$.
   As the quotient is abelian
   (and thus amenable), $\derived{G}\in\calK_\alpha$ implies
   $G\in\calK_{\alpha+1}$, and the result follows inductively.
\end{proof}

\begin{lemma}
   \label{lemma:bootstrap}
   Let $G$ be an extension of $H\in\calK_0$ with cyclic kernel. Then
   $G\in\calK_2$.
\end{lemma}

\begin{proof}
   Let $A$ be the cyclic group, so that the sequence
   \[
     1\longrightarrow A\longrightarrow G
     \stackrel{\pi}{\longrightarrow} H\longrightarrow 1
   \]
   is exact.

   If $H$ is amenable, then $G$ is an amenable group,
   since $A$ is amenable and the class of amenable groups is closed
   under taking extensions.

   Suppose instead that $H$ is free. As free groups have trivial
   second cohomology, $G$ must be the semi-direct product
   $A\rtimes_\phi H$ with $\phi: H\to\Aut A$.

   As before, denote the derived groups of $G$ and $H$ by
   $\derived{G}$ and $\derived{H}$ respectively. Then the following
   diagram commutes with exact rows, and where the vertical homomorphisms
   are just inclusions.
   \[
     \begin{xy}
       \xymatrix{
         1 \ar[r] &
         A\cap\derived{G} \ar[d] \ar[r] &
         \derived{G} \ar[d] \ar[r]^{\pi\vert_{\derived{G}}} &
         \derived{H} \ar[d] \ar[r] &
         1
         \\
         1 \ar[r] &
         A \ar[r] &
         G \ar[r]^\pi &
         H \ar[r] &
         1
       }
     \end{xy}
   \]
   $\derived{G}$ is a semidirect product
   $A\cap\derived{G}\rtimes_{\phi\vert_{\derived{H}}}\derived{H}$,
   and $\phi\vert_{\derived{H}}$ must have image in $\derived{(\Aut   
A)}$.

   As $A$ is cyclic, it has abelian automorphism group,
   and so
   $\derived{(\Aut{A})}$ is trivial. This in turn implies that
   $\phi\vert_{\derived{H}}$ is trivial, and that
   $\derived{G}\cong (A\cap\derived{G})\times\derived{H}$.
   It follows then that $\derived[2]{G}\cong\derived[2]{H}$,
   which is free and thus in $\calK_0$.
   Therefore $G$ is in $\calK_2$.
\end{proof}

\begin{proof}[Proof of Proposition \ref{prop:cyclic}]
   Proceeding by transfinite induction, suppose that for a given
   ordinal $\beta$, $H\in\calK_\alpha$ implies that every extension of
   $H$ with cyclic kernel is in $\calK$ for all $\alpha<\beta$.

   Let $H\in\calK_\beta$, and $G$ an extension of $H$
   with
   \[
     1\longrightarrow A\longrightarrow G
     \stackrel{\pi}{\longrightarrow} H\longrightarrow 1
   \]
   exact and $A$ cyclic.  Choose $\alpha<\beta$ such that
   $H\in\calK_{\alpha+1}$; when $\beta$ is not a
   limit ordinal, one can just let $\alpha=\beta-1$.  There are two
   possibilities: either the group $H$ is a directed union of groups
   $H_j$ in $\calK_\alpha$, $j\in J$; or $H$ is an extension of a group
   $\bar{H}$ in $\calK_\alpha$.

   In the first case, $G=\cup_{j\in J} G_j$ with $G_j=\pi^{-1}(H_j)$.
   Each $G_j$ is an extension of $H_j\in\calK_\alpha$ with cyclic
   kernel, and so $G_j\in\calK$ by the induction hypothesis.
   Therefore $G\in\calK$.

   In the second case, let $B=H/\bar{H}$ be the quotient of the
   extension, and denote the surjection $H\to B$ by $\eta$.
   Let $\bar{G}=\pi^{-1}(\bar{H})$.
   Then $\bar{G}=\ker \eta\pi$ and we get the following commutative
   diagram with exact rows.
   \[
     \begin{xy}
       \xymatrix{
         1 \ar[r] & A \ar@{=}[d] \ar[r] &
         \bar{G} \ar[d] \ar[r]^{\pi\vert_{\bar{G}}} &
         \bar{H} \ar[d] \ar[r] & 1
         \\
         1 \ar[r] & A \ar[r] &
         G \ar[d]^{\eta\pi} \ar[r]^\pi &
         H \ar[d]^\eta \ar[r] & 1
         \\
         & & B \ar@{=}[r] &
         B &
       }
     \end{xy}
   \]
   $G$ is therefore an extension of $\bar{G}$ with amenable quotient   
   $B$,
   while $\bar{G}$ is an extension of $\bar{H}\in\calK_\alpha$ with
   cyclic kernel $A$.
   So $\bar{G}\in\calK$ by the induction hypothesis, and
   so $G\in\calK$.

   The case for $\beta=1$ holds by virtue of
   Lemma~\ref{lemma:bootstrap},
   and so by induction, $H\in\calK$ implies that any extension with
   cyclic kernel of $H$ is in $\calK$.
\end{proof}

\begin{remark}\label{linnell}
   Linnell's class $\calC$ is also closed under taking extensions
   with cyclic kernel,
   by the same argument. $\calC$ can be written as the union of
   classes $\calC_\alpha$ for ordinals $\alpha$:
   \begin{itemize}
   \item
     $\calC_0$ is the class of free groups.
   \item
     $\calC_{\alpha+1}$ is the class of groups which are extensions
     of groups in $\calC_\alpha$ with elementary amenable quotient,
     or are directed unions of groups in $\calC_\alpha$.
   \item
     $\calC_\beta=\bigcup_{\alpha<\beta} \calC_\alpha$ when $\beta$
     is a limit ordinal.
   \end{itemize}
   As the abelian groups are all elementary amenable, the above   
   procedure for $\calK$ can be applied directly to $\calC$.
\end{remark}

\vspace*{10ex}


\begin{thebibliography}{9999}


\bibitem{AS} L.~Ahlfors, L.~Sario,
   Riemann surfaces,
   Princeton Mathematical Series, No. 26
   Princeton University Press, Princeton, N.J. 1960.

\bibitem{At} M.~Atiyah,
   {\em {Elliptic operators, discrete groups and Von
         Neumann algebras}}, Ast\'erisque no. 32-33 (1976) 43-72.

\bibitem{Bel} J.~Bellissard,
   Gap Labeling Theorems for Schr\"{o}dinger's Operators,
   {}From number theory to
   physics (Les Houches, 1989),  538--630, Springer, Berlin, 1992.


\bibitem{Boca} F.~Boca,
   Rotation $C\sp {*}$-algebras and almost Mathieu operators.
   Theta Series in Advanced Mathematics, 1. The Theta Foundation,
   Bucharest, 2001.
   
 \bibitem{Brown} K.~Brown,  
   Cohomology of groups. 
   Graduate Texts in Mathematics, 87. Springer-Verlag, New York-Berlin, 1982.

\bibitem{CHMM98} A.~Carey, K.~Hannabuss, V.~Mathai, P.~McCann,
   Quantum Hall effect on the hyperbolic plane.
   {\em Comm. Math. Phys.} {\bf 190} (1998), no. 3, 629--673.

\bibitem{CHM} A.~Carey, K.~Hannabuss and V.~Mathai,
   Quantum Hall Effect on the hyperbolic plane in the presence
   of disorder,
   {\em Lett. Math. Phys.}, {\bf 47} (1999)
   215--236.


\bibitem{Do} J.~Dodziuk,
   {\em {De Rham-Hodge theory for L$^{2}$-cohomology of
         infinite coverings}},
   Topology {\bf {16}} (1977) 157-165.

\bibitem{dlmsy} J.~Dodziuk, P.~Linnell, V.~Mathai, T.~Schick
   and S.~Yates, Approximating $L^2$-invariants, and the Atiyah
   conjecture,  {\em Communications in Pure and Applied Mathematics},
   {\bf 56}, no. 7, (2003), 839-873.

\bibitem{DM} J.~Dodziuk and V.~Mathai,
   Approximating $L^2$ invariants of amenable covering
   spaces:
   a combinatorial approach, {\em Jour. Func. Anal.}
   {\bf 154} No. 2 (1998) 359--378.

\bibitem{Elek} G.~Elek,
   On the analytic zero divisor conjecture of Linnell,
   {\em Bull. London Math. Soc.} {\bf 35}, no. 2 (2003),  236--238.

\bibitem{G83}
   R.~Grigorchuk,
   On the Milnor problem of group growth. (Russian)
   {\em Dokl. Akad. Nauk SSSR} {\bf 271} (1983), no. 1, 30--33.

\bibitem{GZ}
   R.~Grigorchuk and A.~\.Zuk,
   The lamplighter group as a group generated by a $2$-state automaton
   and its spectrum,
   {\em Geom. Dedicata},  {\bf 87}, (2001),  209--244.

\bibitem{Lin}
   P.~A.~Linnell,
   Division rings and group von {N}eumann algebras,
   {\em Forum Math.}, {\bf 5} (1993) 561--576.

\bibitem{Lu}
   W.~L\"uck,
   Approximating $L^2$ invariants by their finite
   dimensional
   analogues,
   {\em Geom. and Func. Anal.}, {\bf 4} (1994)
   455--481.

\bibitem{MM99} M.~Marcolli, V.~Mathai,
   Twisted index theory on good orbifolds, I:
   noncommutative Bloch theory,
   {\em Communications in Contemporary Mathematics},
   {\bf 1} no. 4 (1999)
   553--587.

\bibitem{MM01} M.~Marcolli, V.~Mathai,
   Twisted index theory on good orbifolds, II: fractional quantum  
   numbers,
   {\em Comm. Math. Phys.}  {\bf 217}, no.1 (2001) 55-87.

\bibitem{MY} V.~Mathai and S.~Yates,
   Approximating spectral invariants
   of Harper operators on graphs, \textit{J. Functional Analysis,}
   {\bf 188}, no. 1 (2002) 111-136.

\bibitem{MSY} V.~Mathai, T.~Schick and S.~Yates,
   Approximating spectral invariants
   of Harper operators on graphs II, \textit{Proc. Amer. Math. Soc.,}
   {\bf 131}, no. 6 (2003), 1917-1923.

\bibitem{Schick(1998a)}
   Thomas Schick,
   \newblock ${L}\sp 2$-determinant class and approximation
   of ${L}\sp 2$-{B}etti numbers.
   \newblock {\em Trans. Amer. Math. Soc.},
   {\bf 353}, no. 8 (2001) 3247--3265.

\bibitem{Sh} M.~Shubin,
   Discrete Magnetic Schr\"odinger operators, {\em Commun. Math. Phys.},
   {\bf 164} (1994), no.2,
   259--275.

\bibitem{Sh2} M.~Shubin,
   von {N}eumann algebras and ${L}\sp 2$ techniques in geometry and topology,
   {\em book in preparation}.

\bibitem{Sun} T.~Sunada,
   A discrete analogue of periodic magnetic Schr\"odinger operators,
   \textit{Contemp. Math.} \textbf{173} (1994), 283--299.



\end{thebibliography}

\end{document}